# On automorphisms of certain free nilpotent-by-abelian Lie algebras

C.E. Kofinas and A.I. Papistas

C.E. Kofinas, Department of Mathematics, University of the Aegean, Karlovassi, 832 00 Samos, Greece. *e-mail:* kkofinas@aegean.gr

A.I. Papistas, Department of Mathematics, Aristotle University of Thessaloniki, 541 24 Thessaloniki, Greece. *e-mail:* apapist@math.auth.gr

**Abstract**

For a positive integer $n$, with $n \geq 4$, let $R_n$ be a free (nilpotent of class 2)-by-abelian and abelian-by-(nilpotent of class 2) Lie algebra of rank $n$. We show that the subgroup of $\mathrm{Aut}(R_n)$ generated by the tame automorphisms and a countably infinite set of explicitly given automorphisms of $R_n$ is dense in $\mathrm{Aut}(R_n)$ with respect to the formal power series topology.

## 1 Introduction

Let $K$ be a field of characteristic 0. By a *Lie algebra* we mean a Lie algebra over $K$. We also write $K$-vector (sub)space instead of vector (sub)space over $K$. For a Lie algebra $L$ and a positive integer $c$, let $\gamma_c(L)$ be the $c$-th term of the lower central series of $L$. We write $L' = \gamma_2(L)$, that is, $L'$ denotes the derived algebra of $L$, and $L'' = (L')'$. We use the left-normed convention for Lie commutators. For a non-negative integer $\kappa$, we write $[a, {}_\kappa b]$ for the left-normed Lie commutator of $a$ and $\kappa$ copies of $b$ (this Lie commutator is interpreted as $a$ if $\kappa = 0$). We write $\mathrm{Aut}(L)$ for the automorphism group of $L$. For a positive integer $n$, with $n \geq 2$, let $L_n$ be a free Lie algebra of rank $n$ with a free generating set $\{x_1, \ldots, x_n\}$. Furthermore $L_n = \bigoplus_{c \geq 1} L_n^c$, where $L_n^c$ is the vector subspace of $L_n$ spanned by all Lie commutators of total degree $c$ in $x_1, \dots, x_n$. If $I$ is a proper fully invariant ideal of $L_n$, then $I \subseteq L_n'$. The general linear group $\mathrm{GL}_n(K)$ acts naturally on $L_n^1$ and we can extend this action so that $\mathrm{GL}_n(K)$ becomes a group of algebra automorphisms of $L_n$. If $f$ belongs to the subalgebra of $L_n$ generated by the set $\{x_2, \dots, x_n\}$, then the endomorphism $\tau_f$ of $L_n$ defined by $\tau_f(x_1) = x_1 + f$, $\tau_f(x_i) = x_i$, $i \neq 1$, is clearly an automorphism of $L_n$. By a result of Cohn [9], $\mathrm{Aut}(L_n)$ is generated by $\mathrm{GL}_n(K)$ and the automorphisms $\tau_f$.

Let $\mathfrak{V}$ be a non-trivial variety of Lie algebras and let $L_n(\mathfrak{V}) = L_n/\mathfrak{V}(L_n)$, where $\mathfrak{V}(L_n)$ denotes the fully invariant ideal of $L_n$ corresponding to $\mathfrak{V}$. Thus $L_n(\mathfrak{V})$ is a relatively free Lie algebra of rank $n$, freely generated by the set $\{y_1, \dots, y_n\}$ where $y_i = x_i + \mathfrak{V}(L_n)$ for $i = 1, \dots, n$. Since $K$ is infinite, we may see by a Vandermonde determinant argument (briefly, a V.d.a) that $\mathfrak{V}(L_n) = \bigoplus_{c\geq 1}(L_n^c \cap \mathfrak{V}(L_n))$ (see, for example, [2, Section 4.2, Theorem 4]). Thus we may write $L_n(\mathfrak{V}) = \bigoplus_{c\geq 1} L_n^c(\mathfrak{V})$, where $L_n^c(\mathfrak{V}) = (L_n^c + \mathfrak{V}(L_n))/\mathfrak{V}(L_n) \cong L_n^c/(L_n^c \cap \mathfrak{V}(L_n))$ and $L_n^c(\mathfrak{V})$ is the vector subspace of $L_n(\mathfrak{V})$ spanned by all Lie commutators of total degree $c$ in $y_1, \dots, y_n$. Since $\mathfrak{V}(L_n) \subseteq L_n'$, the $K$-vector space $L_n^1(\mathfrak{V})$ has a basis $\{y_1, \dots, y_n\}$ and $L_n^1(\mathfrak{V})$ is the natural $K\mathrm{GL}_n(K)$-module. In particular, the general linear group $\mathrm{GL}_n(K)$ acts faithfully on $L_n(\mathfrak{V})$ and we may regard $\mathrm{GL}_n(K)$ as a subgroup of $\mathrm{Aut}(L_n(\mathfrak{V}))$. The natural epimorphism from $L_n$ onto $L_n(\mathfrak{V})$ induces a group homomorphism $\alpha_{n,L}$ from $\mathrm{Aut}(L_n)$ into $\mathrm{Aut}(L_n(\mathfrak{V}))$. An element of $\mathrm{Im}\alpha_{n,L}$ is called *tame*; otherwise it is called *non-tame*. We write $T_{n,L}$ (or, briefly, $T_n$ if it is clear in the context) for $\mathrm{Im}\alpha_{n,L}$.

Bryant and Drensky [5, Section 2] have given a topology on the set of endomorphisms $\mathrm{End}(L_n(\mathfrak{V}))$ of $L_n(\mathfrak{V})$ called the formal power series topology on $\mathrm{End}(L_n(\mathfrak{V}))$. Let $M_n = L_n(\mathfrak{A}\mathfrak{A}) = L_n/L_n''$ be the free metabelian Lie algebra of finite rank $n \geq 2$. In their work [5], Bryant and Drensky have studied $\mathrm{Aut}(M_n)$ by using this topology. They have showed that, for $n \geq 4$, $T_n$ is dense in $\mathrm{Aut}(M_n)$ (see also [14]). In [13, Theorem 1.1 (3)], the above result is extended to a free center-by-metabelian Lie algebra $C_n$ of rank $n \geq 4$ by showing that the subgroup of $\mathrm{Aut}(C_n)$ generated by the tame automorphisms and one more non-tame automorphism is dense in $\mathrm{Aut}(C_n)$. Write $R_n = L_n/(\gamma_3(L_n') + (\gamma_3(L_n))')$ for the free (nilpotent of class 2)-by-abelian and abelian-by-(nilpotent of class 2) Lie algebra of rank $n$. The main purpose of this paper is to show that the subgroup of $\mathrm{Aut}(R_n)$, with $n \geq 4$, generated by the tame automorphisms and a countably infinite set of explicitly given automorphisms of $R_n$ is dense in $\mathrm{Aut}(R_n)$. We point out that the variables involving in the aforementioned automorphisms of $R_n$ are independent upon $n$. In all the above cases (with $n \geq 4$), although a subgroup is dense in the automorphism group it is not clear how far it is from the whole automorphism group. We would like to mention that the method of formal power series topology in the study of automorphisms was first applied by Anick [1]. He considered the formal power series topology on the endomorphisms of the polynomial algebra $K[t_1, \ldots, t_n]$ and he proved that the endomorphisms with invertible Jacobian matrix form a closed subset $J$ and the group of tame automorphisms is dense in $J$.

**Theorem 1.1** *For a positive integer $n$, with $n \geq 4$, let $R_n = L_n/(\gamma_3(L'_n) + (\gamma_3(L_n))')$ be the free (nilpotent of class 2)-by-abelian and abelian-by-(nilpotent of class 2) Lie algebra of rank $n$, freely generated by the set $\{y_1, \ldots, y_n\}$ with $y_i = x_i + (\gamma_3(L'_n) + (\gamma_3(L_n))')$, $i = 1, \ldots, n$. Then the subgroup of $\mathrm{Aut}(R_n)$ generated by the tame automorphisms $T_{R_n}$ and $\omega_1, \omega_2, \ldots$ is dense in $\mathrm{Aut}(R_n)$, where, for $\kappa \geq 1$, $\omega_\kappa(y_1) = y_1 + [y_1, y_2, {}_{(\kappa-1)}y_1, [y_3, y_4]]$ and $\omega_\kappa(y_j) = y_j$, $j = 2, \ldots, n$.*

For a positive integer $n$, with $n \geq 4$, let $L_n/\gamma_3(L'_n)$ be the free (nilpotent of class 2)-by-abelian Lie algebra of rank $n$. For $j \in \{1, \ldots, n\}$, let $q_j = x_j + \gamma_3(L'_n)$. For a positive integer $\kappa$, let $\varphi_\kappa$ be the IA-automorphism of $L_n/\gamma_3(L'_n)$ satisfying the conditions $\varphi_\kappa(q_1) = q_1 + [q_1, q_2, {}_{(\kappa-1)}q_1, [q_3, q_4]]$ and $\varphi_\kappa(q_j) = q_j$, $j \geq 2$ (see Lemma 3.1 below). In Section 4 we show that each automorphism $\varphi_\kappa$ is non-tame by using (right) derivatives defined on the free associative algebra of rank $n$. For $\kappa \in \{1, 2, 3\}$, we use the balanced condition as well. By using this method, we show that $\omega_1, \omega_2$ and $\omega_3$ are non-tame automorphisms of $R_n$. The method does not work for $\omega_\kappa$ with $\kappa \geq 4$. For positive integers $n$ and $\kappa$, with $n \geq 4$, let $G_{n,\kappa}$ be the subgroup of $\mathrm{Aut}(R_n)$ generated by the tame automorphisms $T_{R_n}$ and $\omega_1, \ldots, \omega_\kappa$. We conjecture that, for all $\kappa \geq 1$, $G_{n,\kappa}$ is not dense in $\mathrm{Aut}(R_n)$. In case of the free nilpotent of class 2-by-abelian and abelian-by-nilpotent of class 2 group $F_n/(\gamma_3(F'_n)(\gamma_3(F_n))')$ (where $F_n$ is a free group of finite rank $n$, $n \geq 2$) of rank $n$, it follows from a result of the second author [16] that, for $n \geq 4$, $\mathrm{Aut}(F_n/(\gamma_3(F'_n)(\gamma_3(F_n))'))$ is generated by the tame automorphisms and a countably infinite set of explicitly given automorphisms of $F_n/(\gamma_3(F'_n)(\gamma_3(F_n))')$. Also, the variables involving in the above automorphisms are independent upon $n$.

## 2 Preliminaries

### 2.1 Tame automorphisms

Let $\mathfrak{V}$ be a non-trivial variety of Lie algebras, $L = L_n(\mathfrak{V})$, with $n \geq 2$, freely generated by the set $\{y_i = x_i + \mathfrak{V}(L) : i = 1, \dots, n\}$, and for a positive integer $c$, $L^c = L_n^c(\mathfrak{V})$. Each element $u \in L$ is uniquely written in the form $u = \sum_{k \geq 1} u_k$ with $u_k \in L^k$ for all $k$ and $u_k = 0$ for all but finitely many $k$. We say that $u_k$ is the homogeneous component of $u$ and, if $u_k \neq 0$, we say that $u_k$ is homogeneous of degree $k$. Let $L(y_2, \dots, y_n)$ be the subalgebra of $L$ generated by the set $\{y_2, \dots, y_n\}$. If $u \in L(y_2, \dots, y_n)$, we define $\tau_u \in \mathrm{Aut}(L)$ by $\tau_u(y_1) = y_1 + u$ and $\tau_u(y_i) = y_i$ for $i \geq 2$. For a positive integer $c$, let $L^c(y_2, \dots, y_n)$ be the vector subspace of $L^c$ spanned by all Lie commutators of the form $[y_{i_1}, \dots, y_{i_c}]$ with $i_1, \dots, i_c \in \{2, \dots, n\}$. By the description of $\mathrm{Aut}(L_n)$ (see Introduction), each $\tau_u$ is tame and the group $T_L$ of tame automorphisms of $L$ is generated by $\mathrm{GL}_n(K)$ together with the set of elements $\tau_u$. If $u = u_1 + \cdots + u_m$, where $u_j \in L^j(y_2, \dots, y_n)$, $j = 1, \dots, m$, then $\tau_{u_1} \in \mathrm{GL}_n(K)$ and $\tau_u = \tau_{u_1} \cdots \tau_{u_m}$. Thus $\mathrm{T}_L$ is generated by $\mathrm{GL}_n(K)$ together with those $\tau_u$ for which $u$ is homogeneous of degree at least 2. Since $g\tau_u = (g\tau_u g^{-1})g$ for all $g \in \mathrm{GL}_n(K)$, $\mathrm{T}_L$ can be written as a product of subgroups,

$$\mathrm{T}_L = \langle g\tau_u g^{-1} : g \in \mathrm{GL}_n(K), u \in \bigcup_{\kappa \geq 2} L^\kappa(y_2, \dots, y_n) \rangle \mathrm{GL}_n(K).$$

Let $\mathrm{IA}(L)$ be the kernel of the natural group homomorphism from $\mathrm{Aut}(L)$ into $\mathrm{Aut}(L/L')$. The elements of $\mathrm{IA}(L)$ are called IA-automorphisms of $L$. Write $\mathrm{IT}_L = \mathrm{T}_L \cap \mathrm{IA}(L)$. That is, $\mathrm{IT}_L$ consists of all tame elements of $L$ which induce the identity mapping on $L/L'$. By the modular law and since $\mathrm{IT}_L \cap \mathrm{GL}_n(K) = \{1\}$, we get

$$\mathrm{IT}_L = \langle g\tau_u g^{-1} : g \in \mathrm{GL}_n(K), u \in \bigcup_{\kappa \geq 2} L^\kappa(y_2, \dots, y_n) \rangle.$$

### 2.2 Quotient groups

Let $E_L = \mathrm{End}(L)$ be the semigroup of all (Lie algebra) endomorphisms of $L$. We regard $\mathrm{GL}_n(K)$ as a subgroup of $\mathrm{Aut}(L)$ and so $\mathrm{GL}_n(K) \subseteq E_L$. For the basic terminology and constructions, we refer to [5, Section 2]. For each positive integer $k$, with $k \geq 2$, let $\mathrm{I}_k\mathrm{E}_L$ be the set of endomorphisms of $L$ which induce the identity mapping on $L/\gamma_k(L)$ and write $\mathrm{IE}_L = \mathrm{I}_2\mathrm{E}_L$. Each $\mathrm{I}_k\mathrm{E}_L$ is a subsemigroup of $E_L$. For any element $\phi$ of $\mathrm{I}_k\mathrm{E}_L$ let $\nu_{k,L}(\phi) = (f_1, \dots, f_n)$, where $f_i = (\phi(y_i))_k$ is the homogeneous component of $\phi(y_i)$ of degree $k$, $i = 1, \dots, n$. Thus $\phi(y_i) \equiv y_i + f_i \pmod{\gamma_{k+1}(L)}$, $i = 1, \dots, n$, and $\nu_{k,L}(\phi) \in (L^k)^{\oplus n} = \underbrace{L^k \oplus \cdots \oplus L^k}_{n}$. Each $\nu_{k,L} : \mathrm{I}_k\mathrm{E}_L \longrightarrow (L^k)^{\oplus n}$ is an epimorphism of semigroups. Clearly, for $\phi, \psi \in \mathrm{I}_k\mathrm{E}_L$, $\nu_{k,L}(\phi) = \nu_{k,L}(\psi)$ if and only if $\phi$ and $\psi$ induce the same endomorphism on $L/\gamma_k(L)$. We consider the topology on $L$ corresponding to the series $L \supseteq \gamma_2(L) \supseteq \gamma_3(L) \supseteq \cdots$. Furthermore we give on $\mathrm{E}_L$ the topology of the direct product $L \times L \times \cdots \times L$ of $n$ copies of $L$. The completion $\widehat{L}$ of $L$ with respect to the aforementioned series may be identified with the complete (unrestricted) direct sum $\widehat{\bigoplus}_{i\geq 1} L^i$. It has a natural algebra structure such that $L$ is a subalgebra of $\widehat{L}$. Each element of $\widehat{L}$ may be regarded as an infinite formal sum $u = \sum_{i\geq 1} u_i$ with $u_i \in L^i$ for all $i$. For each $k \geq 1$, let $\widehat{\gamma}_k(L)$ be the set of all such elements $u$ with $u_i = 0$ for $i < k$. That is, $\widehat{\gamma}_k(L)$ is the completion of $\gamma_k(L)$. If $w_1, \dots, w_n$ are arbitrary elements of $\widehat{L}$, then there is a unique continuous endomorphism $\phi$ of $\widehat{L}$ such that $\phi(y_i) = w_i$, $i = 1, \dots, n$. Let $\widehat{E}_L$ be the semigroup of all continuous endomorphisms of $\widehat{L}$. Since each element $\phi$ of $\widehat{E}_L$ corresponds uniquely to an element $(\phi(y_1), \dots, \phi(y_n))$ of the direct product $\widehat{L} \times \cdots \times \widehat{L}$ of $n$ copies of $\widehat{L}$, the set $\widehat{E}_L$ with the topology of this direct product may be identified with the completion of $\mathrm{E}_L$.

For $k \geq 2$, let $\mathrm{I}_k\widehat{\mathrm{E}}_L$ be the set of all elements of $\widehat{E}_L$ which induce the identity map on $\widehat{L}/\widehat{\gamma}_k(L)$. Thus $E_L \cap \mathrm{I}_k\widehat{\mathrm{E}}_L = \mathrm{I}_k\mathrm{E}_L$. We also write $\mathrm{I}\widehat{\mathrm{E}}_L = \mathrm{I}_2\widehat{\mathrm{E}}_L$. As shown in [5, Lemma 2.2], $\mathrm{I}\widehat{\mathrm{E}}_L$ is a group. Thus each $\mathrm{I}_k\widehat{\mathrm{E}}_L$ is a normal subgroup of $\mathrm{I}\widehat{\mathrm{E}}_L$. For each $k$, we extend the semigroup homomorphism $\nu_{k,L} : \mathrm{I}_k\mathrm{E}_L \longrightarrow (L^k)^{\oplus n}$ to a group homomorphism $\nu_{k,L} : \mathrm{I}_k\widehat{\mathrm{E}}_L \longrightarrow (L^k)^{\oplus n}$ in the obvious way. Thus $\nu_{k,L}$ induces a group isomorphism $\overline{\nu}_{k,L} : \mathrm{I}_k\widehat{\mathrm{E}}_L/\mathrm{I}_{k+1}\widehat{\mathrm{E}}_L \longrightarrow (L^k)^{\oplus n}$. For each $k$, we write $\bar{\mathrm{I}}_k\mathrm{E}_L = \mathrm{I}_k\widehat{\mathrm{E}}_L/\mathrm{I}_{k+1}\widehat{\mathrm{E}}_L = (\mathrm{I}_k\mathrm{E}_L)(\mathrm{I}_{k+1}\widehat{\mathrm{E}}_L)/\mathrm{I}_{k+1}\widehat{\mathrm{E}}_L$. Since $(L^k)^{\oplus n}$ is a vector space, by using the map $\overline{\nu}_{k,L}$, we give $\bar{\mathrm{I}}_k\mathrm{E}_L$ the structure of a vector space so that $\overline{\nu}_{k,L}$ is a vector space isomorphism. More explicitly, if $\overline{\phi} \in \bar{\mathrm{I}}_k\mathrm{E}_L$ is represented by $\phi \in \mathrm{I}_k\mathrm{E}_L$ and if $a \in K$ then $a\overline{\phi}$ is represented by the endomorphism $\phi_1$ defined by $\phi_1(y_i) = y_i + af_i$ for all $i$, where $\nu_{k,L}(\phi) = (f_1, \dots, f_n)$. The vector space direct sum $\mathcal{L}(E_L) = \bigoplus_{k\geq 2} \bar{\mathrm{I}}_k\mathrm{E}_L$ has the structure of a graded Lie algebra with $\bar{\mathrm{I}}_k\mathrm{E}_L$ as component of degree $k-1$ in the grading and Lie multiplication given by $[\phi\mathrm{I}_{j+1}\widehat{\mathrm{E}}_L, \psi\mathrm{I}_{k+1}\widehat{\mathrm{E}}_L] = (\phi^{-1}\psi^{-1}\phi\psi)\phi\mathrm{I}_{j+k}\widehat{\mathrm{E}}_L$ for all $\phi \in \mathrm{I}_j\widehat{\mathrm{E}}_L$, $\psi \in \mathrm{I}_k\widehat{\mathrm{E}}_L$, $j, k \geq 2$. The group $\mathrm{GL}_n(K)$ acts on $\mathcal{L}(E_L)$ as a group of Lie algebra automorphisms (see [5, Proposition 2.5]).

As noticed before, $L^1$ is the natural $K\mathrm{GL}_n(K)$-module with basis $\{y_1, \dots, y_n\}$. It is convenient to regard elements of $\mathrm{GL}_n(K)$ as $n \times n$ matrices corresponding to the ordered basis $\{y_1, \ldots, y_n\}$ of $L^1$. Since $\mathrm{GL}_n(K) \subseteq \mathrm{E}_L \subseteq \widehat{\mathrm{E}}_L$, $\mathrm{GL}_n(K)$ acts by conjugation on $\widehat{\mathrm{E}}_L$. It is easily verified that each $\mathrm{I}_k\mathrm{E}_L$ is $\mathrm{GL}_n(K)$-invariant and that if $\phi$ and $\psi$ are elements of $\mathrm{I}_k\mathrm{E}_L$ satisfying $\nu_{k,L}(\phi) = \nu_{k,L}(\psi)$ then $\nu_{k,L}(g\phi g^{-1}) = \nu_{k,L}(g\psi g^{-1})$ for all $g \in \mathrm{GL}_n(K)$. Thus $\mathrm{GL}_n(K)$ acts on $\bar{\mathrm{I}}_k\mathrm{E}_L$. More precisely, $g * \overline{\phi} = \overline{g\phi g^{-1}}$ for all $g \in \mathrm{GL}_n(K)$. It is also easy to see that the action of $\mathrm{GL}_n(K)$ commutes with multiplication by elements of $K$. Thus $\bar{\mathrm{I}}_k\mathrm{E}_L$ is a $K\mathrm{GL}_n(K)$-module. The action of $\mathrm{GL}_n(K)$ on $\bar{\mathrm{I}}_k\mathrm{E}_L$ is most easily written down using the map $\nu_{k,L}$. Let $\phi \in \mathrm{I}_k\mathrm{E}_L$, $g \in \mathrm{GL}_n(K)$ and $\nu_{k,L}(\phi) = (f_1, \ldots, f_n)$. Then $\nu_{k,L}$ maps $g\phi g^{-1}$ to $(gf_1, \ldots, gf_n)g^{-1}$. Here $gf_i$ is calculated in the $\mathrm{GL}_n(K)$-module $L^k$, $g^{-1}$ is regarded an $n \times n$ matrix, and multiplication by $g^{-1}$ is multiplication of a $1 \times n$ matrix by an $n \times n$ matrix.

For any subgroup $G$ of $\mathrm{Aut}(L)$, we write $\mathrm{I}_k\mathrm{G} = G \cap \mathrm{I}_k\widehat{\mathrm{E}}_L$, $k \geq 2$, and $\mathrm{IG} = \mathrm{I}_2\mathrm{G}$. Thus $\mathrm{I}_k\mathrm{G}$ is the set of elements of $G$ which induce the identity map on $L/\gamma_k(L)$ and is a normal subgroup of $G$. We also write $\bar{\mathrm{I}}_k\mathrm{G} = \mathrm{I}_k\mathrm{G}(\mathrm{I}_{k+1}\widehat{\mathrm{E}}_L)/\mathrm{I}_{k+1}\widehat{\mathrm{E}}_L$. Notice that if $G_1$ and $G_2$ are subgroups of $\mathrm{Aut}(L)$ with $G_1 \subseteq G_2$, we have $\bar{\mathrm{I}}_k\mathrm{G}_1 \subseteq \bar{\mathrm{I}}_k\mathrm{G}_2$. Let $G$ be a subgroup of $\mathrm{Aut}(L)$ such that $\mathrm{GL}_n(K) \subseteq G$. It has been shown in [5, Proposition 2.6] that, for $k \geq 2$, $\bar{\mathrm{I}}_k\mathrm{G}$ is a $K\mathrm{GL}_n(K)$-submodule of $\bar{\mathrm{I}}_k\mathrm{E}_L$, and $\mathcal{L}(G) = \bigoplus_{k\geq 2} \bar{\mathrm{I}}_k G$ is a graded Lie algebra which is a $\mathrm{GL}_n(K)$-invariant subalgebra of $\mathcal{L}(E_L)$ (see [5, Proposition 2.7]). Let $G_1$ and $G_2$ be subgroups of $\mathrm{Aut}(L)$ such that $\mathrm{GL}_n(K) \subseteq G_1 \subseteq G_2$. It follows from [5, Corollary 2.9] that $G_1$ is dense in $G_2$ if and only if $\bar{\mathrm{I}}_k G_1 = \bar{\mathrm{I}}_k G_2$ for all $k \geq 2$.

### 2.3 $K\mathrm{GL}_n(K)$-modules

We summarize some information about $K\mathrm{GL}_n(K)$-modules, particularly finite dimensional polynomial $K\mathrm{GL}_n(K)$-modules. For basic facts and definitions for polynomial modules we refer to [12]. A finite direct sum, a finite tensor product, submodules and quotient modules of polynomial $K\mathrm{GL}_n(K)$-modules are polynomial ones. Every polynomial $K\mathrm{GL}_n(K)$-module is a direct sum of homogeneous polynomial $K\mathrm{GL}_n(K)$-submodules and so any irreducible polynomial $K\mathrm{GL}_n(K)$-module is homogeneous. Every homogeneous polynomial $K\mathrm{GL}_n(K)$-module is completely reducible. Every polynomial $K\mathrm{GL}_n(K)$-module is a direct sum of irreducible ones. The irreducible polynomial $K\mathrm{GL}_n(K)$-modules are indexed (up to isomorphism) by the $n$-tuples of non-negative integers $\lambda = (\lambda_1, \ldots, \lambda_n)$ with $\lambda_1 \geq \cdots \geq \lambda_n$. For $\lambda = (\lambda_1, \ldots, \lambda_n)$, with $\lambda_1 \geq \lambda_2 \geq \cdots \geq \lambda_n \geq 0$, the irreducible polynomial $K\mathrm{GL}_n(K)$-module corresponding to $\lambda$ will be denoted by $[\lambda]$ or $[\lambda_1, \ldots, \lambda_n]$. Let $k$ be a positive integer. An $n$-tuple $\lambda = (\lambda_1, \ldots, \lambda_n)$ with $\lambda_1 \geq \cdots \geq \lambda_n$ and $\lambda_1 + \cdots + \lambda_n = k$ is called a partition of $k$ into $n$ parts. When writing partitions we shall make use of standard abbreviations. Thus, for example, $(3, 2, 2, 1, 1, 1, 0)$ is written $(3, 2^2, 1^3)$.

Let $\mathfrak{V}$ be a non-trivial variety of Lie algebras and $n, c$ be positive integers, with $n, c \geq 2$. It is easily verified that $L_n^c(\mathfrak{V})$ is a homogeneous polynomial $K\mathrm{GL}_n(K)$-module of degree $c$. For $g \in \mathrm{GL}_n(K)$ and $(u_1, \ldots, u_n) \in L_n^c(\mathfrak{V})^{\oplus n}$, we define $g(u_1, \ldots, u_n) = (gu_1, \ldots, gu_n)g^{-1}$, where $gu_i$ is calculated in the polynomial $K\mathrm{GL}_n(K)$-module $L_n^c(\mathfrak{V})$, $g^{-1}$ is regarded as an $n \times n$ matrix and multiplication by $g^{-1}$ is multiplication of a $1 \times n$ matrix by an $n \times n$ matrix. It is easily verified that $\mathrm{GL}_n(K)$ acts on $L_n^c(\mathfrak{V})^{\oplus n}$. The action of $\mathrm{GL}_n(K)$ on $L_n^c(\mathfrak{V})$ commutes with multiplication by elements of $K$ and so $L_n^c(\mathfrak{V})^{\oplus n}$ is a $K\mathrm{GL}_n(K)$-nodule. Clearly $L_n^c(\mathfrak{V})^{\oplus n}$ is isomorphic to $\bar{\mathrm{I}}_c\mathrm{E}_L$ as $K\mathrm{GL}_n(K)$-module for all $n, c \geq 2$. It follows from [5, Theorem 1.1, Corollary 1.2] that, for all $n, c \geq 2$,

$$L_n^c(\mathfrak{V})^{\oplus n} \cong (\det)^{-1} \otimes_K [1^{n-1}] \otimes_K L_n^c(\mathfrak{V})$$

as $K\mathrm{GL}_n(K)$-modules, where $(\det)^{-1}$ is a one-dimensional $K\mathrm{GL}_n(K)$-module which affords the representation $g \longrightarrow (\det)^{-1}$ for all $g \in \mathrm{GL}_n(K)$ where $\det g$ is the determinant of $g$.

Let $V$ be a homogeneous polynomial $K\mathrm{GL}_n(K)$-module of degree $k$. For $g \in \mathrm{GL}_n(K)$ and $(v_1, \ldots, v_n) \in V^{\oplus n}$, we define $g(v_1, \ldots, v_n) = (gv_1, \ldots, gv_n)g^{-1}$, where $gv_i$ is calculated in the module $V$, $g^{-1}$ is regarded as an $n \times n$ matrix and multiplication by $g^{-1}$ is multiplication of a $1 \times n$ matrix by an $n \times n$ matrix. The group $\mathrm{GL}_n(K)$ acts on $V^{\oplus n}$. The action of $\mathrm{GL}_n(K)$ on $V$ commutes with multiplication by elements of $K$ and so $V^{\oplus n}$ is a $K\mathrm{GL}_n(K)$-nodule. Let $[1]^*$ be the vector space of $1 \times n$ row-vectors over $K$ regarded as a left $K\mathrm{GL}_n(K)$-module in which, for each $g \in \mathrm{GL}_n(K)$, $g$ acts as right multiplication by $g^{-1}$ and regard $V \otimes_K [1]^*$ as a $K\mathrm{GL}_n(K)$-module under the diagonal action of $\mathrm{GL}_n(K)$. As observed in [5, at the bottom of page 1138], the map $(v_1, \ldots, v_n) \mapsto v_1 \otimes (1, 0, \ldots, 0) + \cdots v_n \otimes (0, \ldots, 0, 1)$ determines a $K\mathrm{GL}_n(K)$-module isomorphism from $V^{\oplus n}$ to $V \otimes_K [1]^*$. Since $[1]^* \cong (\det)^{-1} \otimes [1^{n-1}]$, we have

$$V^{\oplus n} \cong (\det)^{-1} \otimes_K [1^{n-1}] \otimes_K V. \tag{2.1}$$

## 3 Certain nilpotent-by-abelian Lie algebras

Throughout this paper, for $n \geq 4$, we write $I_n = \gamma_3(L_n') + (\gamma_3(L_n))'$ and $R_n = L_n/I_n$. Set $y_i = x_i + I_n$ with $i = 1, \ldots, n$. Thus $\{y_1, \ldots, y_n\}$ is a free generating set of $R_n$. We order this generating set as $y_1 < y_2 < \cdots < y_n$. For a positive integer $d$, let $R_n^d = (L_n^d + I_n)/I_n$. Since $I_n = \bigoplus_{d \geq 1}(I_n \cap L_n^d)$, we have $R_n = \bigoplus_{d \geq 1} R_n^d$.

### 3.1 Analysis of $R_n$

For any Lie algebra $L$ and each $v \in L$, $\mathrm{ad}v : L \longrightarrow L$ is defined by $u(\mathrm{ad}v) = [u, v]$ for all $u \in L$. Let $D_n = K[t_1, \ldots, t_n]$ be the (commutative, associative, unitary) polynomial algebra over $K$ freely generated by the variables $t_1, \ldots, t_n$. For $k \geq 0$, we write $D_n^k$ for the homogeneous component of $D_n$ of degree $k$ and $D_n = \bigoplus_{k \geq 0} D_n^k$, with $D_n^0 = K$. Working modulo $R_n'' = L_n''/I_n$, every element of the derived algebra $R_n' = L_n'/I_n$ of $R_n$ may be written in the form $\sum_{1 \leq i,j \leq n}[y_i, y_j]f_{ij}(\mathrm{ad}y_1, \ldots, \mathrm{ad}y_n)$, where $f_{ij}(t_1, \ldots, t_n) \in D_n$ for all $i, j$. Writing $[y_i, y_j]f_{ij}(\mathrm{ad}y_1, \ldots, \mathrm{ad}y_n) = [y_i, y_j] \cdot f_{ij}(t_1, \ldots, t_n)$ for all $i, j$, we have

$$\sum_{1 \leq i,j \leq n} [y_i, y_j]f_{ij}(\mathrm{ad}y_1, \ldots, \mathrm{ad}y_n) = \sum_{1 \leq i,j \leq n} [y_i, y_j] \cdot f_{ij}(t_1, \ldots, t_n).$$

By using the Jacobi identity in the form $[x, y, z] = [x, z, y] + [x, [y, z]]$ and since $\gamma_3(R_n') = (\gamma_3(R_n))' = \{0\}$, we have $R_n''$ is spanned by all elements of the form $[[y_{i_1}, y_{i_2}, \ldots, y_{i_\kappa}], [y_{j_1}, y_{j_2}]]$, where $\kappa \geq 2$, $i_1 > i_2 \leq i_3 \leq \cdots \leq i_\kappa$ and $j_1 > j_2$. Furthermore $R_n'' = \bigoplus_{d \geq 4}(R_n'')^d$, where $(R_n'')^d$ is spanned by the set $\{[y_{i_1}, \ldots, y_{i_{d-2}}, [y_{j_1}, y_{j_2}]] : i_1 > i_2 \leq i_3 \leq \cdots \leq i_{d-2}, j_1 > j_2\}$. For $i \in \{1, \ldots, n\}$, let $z_i = x_i + L_n''$. Thus the set $\{z_1, \ldots, z_n\}$ is a free generating set of $M_n$. For a positive integer $d$, let $M_n^d$ be the $K$-vector subspace of $M_n$ spanned by all Lie commutators $[z_{i_1}, \ldots, z_{i_d}]$. Thus $M_n = \bigoplus_{d \geq 1} M_n^d$.

For each $d \geq 1$, the natural surjection from $R_n$ onto $M_n$ induces a $K\mathrm{GL}_n(K)$-module epimorphism $\mu_{d,R_n} : R_n^d \to M_n^d$. In particular, $\mu_{d,R_n}(u + I_n) = u + L_n''$ for all $u \in L_n^d$. We point out that the kernel of $\mu_{d,R_n}$ is $(R_n'')^d$. Since $R_n^d$ is completely reducible, there exists a polynomial $K\mathrm{GL}_n(K)$-submodule $\Omega_n^d$ of $R_n^d$ such that $R_n^d = \Omega_n^d \oplus (R_n'')^d$ and $\Omega_n^d \cong M_n^d$ as $K\mathrm{GL}_n(K)$-module. Notice that $R_n^d = \Omega_n^d$ for $d \in \{1, 2, 3\}$. It is well known that $M_n^1 \cong [1]$ and, for $d \geq 2$, $M_n^d \cong [d-1, 1]$. For $n \geq 4$, $(R_n'')^4 \cong [2, 1^2]$ (see [19]) and, if $d \geq 5$, then $(R_n'')^d \cong [d-3, 1] \otimes [1^2]$. The tensor product of polynomial modules can be calculated by means of the Littlewood-Richardson rule (see [15]). Thus, for $n \geq 4$ and $d \geq 5$,

$$(R''_n)^d \cong [d-2,2] \oplus [d-2,1^2] \oplus [d-3,2,1] \oplus [d-3,1^3]. \tag{3.1}$$

### 3.2 Automorphisms of $R_n$

The proof of the following result is based on the arguments given in [7, Lemma 3.1] (see, also, [11]).

**Lemma 3.1** *Let $V$ be a fully invariant ideal of $L_n$, with $n \geq 2$, and $V \subseteq L''_n$ such that $L'_n/V$ is nilpotent. Let $\lambda$ be an endomorphism of $L_n$. Then $\lambda$ induces an automorphism of $L_n/V$ if and only if $\lambda$ induces an automorphism of $L_n/L''_n$.*

The natural mapping from $R_n$ onto $M_n$ induces a group homomorphism $\pi_n$ from $\mathrm{Aut}(R_n)$ into $\mathrm{Aut}(M_n)$. By Lemma 3.1, $\pi_n$ is surjective. By a result of Bahturin and Nabiyev [3], $M_n$ has non-tame automorphisms for all $n \geq 2$ and so there exist non-tame automorphisms of $R_n$ for all $n \geq 2$. Write $\mathrm{Ker}\pi_n$ for the kernel of $\pi_n$. Since $R_n/R'_n \cong M_n/M'_n \cong L_n/L'_n$ as $K$-vector spaces, the restriction $\widetilde{\pi}_n$ of $\pi_n$ on $\mathrm{IA}(R_n)$ is an epimorphism from $\mathrm{IA}(R_n)$ onto $\mathrm{IA}(M_n)$. Since $\mathrm{Ker}\widetilde{\pi}_n = \mathrm{Ker}\pi_n \subseteq \mathrm{IA}(R_n)$, we have $\mathrm{IA}(R_n)/\mathrm{Ker}\pi \cong \mathrm{IA}(M_n)$ as groups. Clearly an element $\phi \in \mathrm{Ker}\pi_n$ satisfies the conditions $\phi(y_i) = y_i + u_i$, $u_i \in R''_n$, $i = 1, \ldots, n$. Since $\gamma_3(R'_n) = \{0\}$, the elements of $\mathrm{Ker}\pi_n$ act trivially on $R''_n$, and this fact will be used repeatedly below without further comment. Furthermore, since $R''_n$ is an abelian Lie algebra, $\mathrm{Ker}\pi_n$ is an abelian group. For a positive integer $k \geq 4$, let $\Pi_{n,k}$ be the subset of $\mathrm{Ker}\pi_n$ consisting of all IA-automorphisms $\phi$ satisfying the conditions $\phi(y_i) = y_i + u_i$, $u_i \in (R''_n)^k$, $i = 1, \ldots, n$. Since $R''_n = \bigoplus_{k \geq 4} (R''_n)^k$ and $\mathrm{Ker}\pi_n$ is an abelian group, each $\Pi_{n,k}$ is a subgroup of $\mathrm{Ker}\pi_n$ and $\mathrm{Ker}\pi_n = \bigoplus_{k \geq 4} \Pi_{n,k}$.

For $k \geq 4$, the mapping $\eta_k : \Pi_{n,k} \to [(R''_n)^k]^{\oplus n}$ defined by $\eta_k(\phi) = (u_1, \ldots, u_n)$, where $\phi(y_i) = y_i + u_i$, $u_i \in (R''_n)^k$, $i \in \{1, \ldots, n\}$, is clearly an isomorphism of abelian groups. By using the group isomorphism $\eta_k$, we give each $\Pi_{n,k}$ a structure of a $K$-vector space so that $\eta_k$ is a vector space isomorphism. Notice that if $\phi \in \Pi_{n,k}$ and $a \in K$, then $a \cdot \phi$ is defined by $(a \cdot \phi)(y_i) = y_i + au_i$ for all $i \in \{1, \ldots, n\}$. Since $\mathrm{GL}_n(K) \leq \mathrm{Aut}(R_n)$, we can let $\mathrm{GL}_n(K)$ act by conjugation on $\mathrm{IA}(R_n)$. It is easily verified that $g\phi g^{-1} \in \Pi_{n,k}$ for all $g \in \mathrm{GL}_n(K)$ and $\phi \in \Pi_{n,k}$. Thus $\mathrm{GL}_n(K)$ acts on $\Pi_{n,k}$. It is also easy to see that the action of $\mathrm{GL}_n(K)$ on $\Pi_{n,k}$ commutes with multiplication by elements of $K$. Thus each $\Pi_{n,k}$ is a $K\mathrm{GL}_n(K)$-module. The action of $\mathrm{GL}_n(K)$ on each $\Pi_{n,k}$ is most easily written down using the isomorphism $\eta_k$. Let $\phi \in \Pi_{n,k}$, $g \in \mathrm{GL}_n(K)$ and $\eta_k(\phi) = (u_1, \ldots, u_n)$ with $u_1, \ldots, u_n \in (R''_n)^k$. Then $\eta_k(g\phi g^{-1}) = (gu_1, \ldots, gu_n)g^{-1}$. Clearly $\Pi_{n,k}$ and $[(R''_n)^k]^{\oplus n}$ are isomorphic as $K\mathrm{GL}_n(K)$-modules. By (2.1), for all $n, k \geq 4$,

$$\Pi_{n,k} \cong (\det)^{-1} \otimes_K [1^{n-1}] \otimes_K (R''_n)^k \tag{3.2}$$

as $K\mathrm{GL}_n(K)$-modules. Since $\mathrm{Ker}\pi_n = \bigoplus_{k\geq 4} \Pi_{n,k}$ and each $\Pi_{n,k}$ has a strucure of $K\mathrm{GL}_n(K)$-module, we may consider $\mathrm{Ker}\pi_n$ as $K\mathrm{GL}_n(K)$-module.

For $n, k \geq 4$, let $\overline{\nu}_{k,R_n}$ be the $K\mathrm{GL}_n(K)$-module isomorphism from $\bar{\mathrm{I}}_k\mathrm{E}_{R_n}$ to $(R^k_n)^{\oplus n}$ sending $\overline{\phi}_R$ to $(u_1, \ldots, u_n)$, where $\phi_R(y_i) = y_i + u_i$, $u_i \in R^k_n$, $i = 1, \ldots, n$. Let $\overline{\mu}_{k,R_n}$ be the $K\mathrm{GL}_n(K)$-module epimorphism from $(R^k_n)^{\oplus n}$ onto $(M^k_n)^{\oplus n}$ induced by the $K\mathrm{GL}_n(K)$-module epimorphism $\mu_{k,R_n}$ from $R^k_n$ onto $M^k_n$. That is, for all $r_1, \ldots, r_n \in R^k_n$ and $g \in \mathrm{GL}_n(K)$, $\overline{\mu}_{k,R_n}(g(r_1, \ldots, r_n)) = g(\mu_{k,R_n}(r_1), \ldots, \mu_{k,R_n}(r_n))$. Define

$$\overline{\vartheta}_{k,R_n} = (\overline{\nu}_{k,M_n})^{-1}\, \overline{\mu}_{k,R_n}\, \overline{\nu}_{k,R_n}.$$

Since $\overline{\nu}_{k,M_n}$ is a $K\mathrm{GL}_n(K)$-module isomorphism from $\bar{\mathrm{I}}_k\mathrm{E}_{M_n}$ to $(M^k_n)^{\oplus n}$, we get $\overline{\vartheta}_{k,R_n}$ is a $K\mathrm{GL}_n(K)$-module epimorphism from $\bar{\mathrm{I}}_k\mathrm{E}_{R_n}$ onto $\bar{\mathrm{I}}_k\mathrm{E}_{M_n}$. In particular if $\overline{\phi}_R \in \bar{\mathrm{I}}_k\mathrm{E}_{R_n}$, with $\phi_R \in \mathrm{I}_k\mathrm{E}_{R_n}$ satisfying the conditions $\phi_R(y_i) = y_i + u_i$, $u_i \in R^k_n$, $i = 1, \ldots, n$, then $\overline{\vartheta}_{k,R_n}(\overline{\phi}_R) = \overline{\phi}_M$, where $\phi_M \in \mathrm{I}_k\mathrm{E}_{M_n}$ such that $\phi_M(z_i) = z_i + \mu_{k,R_n}(u_i)$, $i = 1, \ldots, n$. Moreover, for all $g \in \mathrm{GL}_n(K)$, $\overline{\vartheta}_{k,R_n}(\overline{g\phi_R g^{-1}}) = \overline{g\phi_M g^{-1}}$.

**Lemma 3.2** *For positive integers $n, k$, with $n, k \geq 4$, let $\mathrm{Ker}\overline{\vartheta}_{k,R_n}$ be the kernel of $\overline{\vartheta}_{k,R_n}$. Then $\mathrm{Ker}\overline{\vartheta}_{k,R_n} \subseteq \bar{\mathrm{I}}_k\mathrm{A}_{R_n}$ and, as $K\mathrm{GL}_n(K)$-module,*

$$\mathrm{Ker}\overline{\vartheta}_{k,R_n} \cong \Pi_{n,k} \cong (\det)^{-1} \otimes_K [1^{n-1}] \otimes_K (R''_n)^k.$$

*Proof.* For $n, k \geq 4$, $\mathrm{Ker}\overline{\vartheta}_{k,R_n}$ consists of those elements $\overline{\phi}_R$ of $\bar{\mathrm{I}}_k\mathrm{E}_{R_n}$ where $\phi_R \in \mathrm{I}_k\mathrm{E}_{R_n}$ satisfying the conditions $\phi_R(y_i) = y_i + v_i$, $v_i \in (R''_n)^k$, $i = 1, \ldots, n$. By Lemma 3.1, each $\phi_R \in \mathrm{Aut}(R_n)$ and so $\mathrm{Ker}\overline{\vartheta}_{k,R_n} = \{\overline{\phi}_R \in \bar{\mathrm{I}}_k\mathrm{A}_{R_n} : \phi_R \in \Pi_{n,k}\} \subseteq \bar{\mathrm{I}}_k\mathrm{A}_{R_n}$. It is clear enough that $\mathrm{Ker}\overline{\vartheta}_{k,R_n}$ is a $K\mathrm{GL}_n(K)$-submodule of $\bar{\mathrm{I}}_k\mathrm{A}_{R_n}$, and it is isomorphic to $\Pi_{n,k}$ as $K\mathrm{GL}_n(K)$-module. By (3.2) we obtain the required result. □

By (3.1), (3.2), Lemma 3.2 and the Littlewood-Richardson rule, we obtain the following result.

**Lemma 3.3** *Let $n$ and $k$ be positive integers, with $n, k \geq 4$. Then, as $K\mathrm{GL}_n(K)$-modules,*

1. $\mathrm{Ker}\overline{\vartheta}_{4,R_n} \cong ((\det)^{-1} \otimes_K [3, 2^2, 1^{n-4}]) \oplus [2, 1] \oplus [1^3]$.

2. $\mathrm{Ker}\overline{\vartheta}_{5,R_4} \cong ((\det)^{-1} \otimes [4, 3, 1]) \oplus ((\det)^{-1} \otimes [4, 2^2]) \oplus ((\det)^{-1} \otimes [3^2, 2]) \oplus 2[3, 1] \oplus 2[2^2] \oplus 3[2, 1^2] \oplus [1^4]$.

3. *For $k \geq 6$,*

$$\begin{aligned} \mathrm{Ker}\overline{\vartheta}_{k,R_4} \cong\ & ((\det)^{-1} \otimes [k-1, 3, 1]) \oplus ((\det)^{-1} \otimes [k-1, 2^2]) \oplus ((\det)^{-1} \otimes [k-2, 3, 2]) \\ & \oplus\ 2[k-2, 1] \oplus 2[k-3, 2] \oplus 3[k-3, 1^2] \oplus [k-4, 2, 1] \oplus [k-4, 1^3]. \end{aligned}$$

4. *For $n, k \geq 5$,*

$$\begin{aligned} \mathrm{Ker}\overline{\vartheta}_{k,R_n} \cong\ & ((\det)^{-1} \otimes [k-1, 3, 1^{n-3}]) \oplus ((\det)^{-1} \otimes [k-1, 2^2, 1^{n-4}]) \oplus \\ & ((\det)^{-1} \otimes [k-2, 3, 2, 1^{n-4}]) \oplus ((\det)^{-1} \otimes [k-2, 2^3, 1^{n-5}]) \oplus \\ & 2[k-2, 1] \oplus 2[k-3, 2] \oplus 3[k-1, 1^2] \oplus [k-4, 2, 1] \oplus [k-4, 1^3]. \end{aligned}$$

### 3.3 A further analysis of $\mathrm{Ker}\overline{\vartheta}_{k,R_n}$ as $K\mathrm{GL}_n(K)$-module

For $i \in \{1, \ldots, n\}$ and $u \in R''_n$, we write $\phi_{i,u}(y_i) = y_i + u$ and $\phi_{i,u}(y_j) = y_j$ for $j \neq i$. By Lemma 3.1, $\phi_{i,u}$ is an automorphism of $R_n$. In particular, $\phi_{i,u} \in \mathrm{Ker}\pi_n$ for all $i \in \{1, \ldots, n\}$ and $u \in R''_n$.

**Lemma 3.4** *For a permutation $\rho$ of $\{1, \ldots, n\}$, let $\tau_\rho$ be the tame automorphism of $R_n$ satisfying the conditions $\tau_\rho(y_i) = y_{\rho(i)}$ for all $i \in \{1, \ldots, n\}$. Then $\tau_\rho \phi_{i,u} \tau_\rho^{-1} = \phi_{\rho(i), \tau_\rho(u)}$ for all $i \in \{1, \ldots, n\}$ and $u \in R''_n$.*

*Proof.* This is straightforward. □

**Lemma 3.5** *For $n, k \geq 4$, $\mathrm{Ker}\overline{\vartheta}_{k,R_n}$ is generated as $K\mathrm{GL}_n(K)$-module by all automorphisms of the form $\overline{\phi}$ with $\phi(y_1) = y_1 + u$ and $\phi(y_i) = y_i$, $i \geq 2$, $u \in (R''_n)^k$.*

*Proof.* Recall that $\mathrm{Ker}\overline{\vartheta}_{k,R_n} = \{\overline{\phi}_R \in \bar{\mathrm{I}}_k\mathrm{A}_{R_n} : \phi_R \in \Pi_{n,k}\}$. Consider the tame automorphisms $\tau_\rho$ of Lemma 3.4 as elements of $\mathrm{GL}_n(K)$. By the action of $\mathrm{GL}_n(K)$ on $\bar{\mathrm{I}}_k\mathrm{E}_{R_n}$ and Lemma 3.4, we obtain the required result. □

For positive integers $i$ and $j$, with $i, j \in \{1, \ldots, n\}$ and $i \neq j$, we write $\tau_{i,j,a}$ and $\sigma_{i,j}$ for the tame automorphisms of $R_n$ satisfying the conditions $\tau_{i,j,a}(y_i) = y_i + a y_j$, $a \in K \backslash \{0\}$, $\tau_{i,j,a}(y_k) = y_k$, for $k \neq i$, and $\sigma_{i,j}(y_i) = y_j$, $\sigma_{i,j}(y_j) = y_i$, $\sigma_{i,j}(y_k) = y_k$ for $k \neq i, j$. If $a = 1$ we write $\tau_{i,j}$ instead of $\tau_{i,j,1}$. The above automorphisms of $R_n$ may be regarded as elements of $\mathrm{GL}_n(K)$ as well in a natural way.

Recall that $D_n = K[t_1, \ldots, t_n] = \bigoplus_{k \geq 0} D_n^k$, with $D_n^0 = K$ and, for $k \geq 1$, $D_n^k$ is the vector subspace of $D_n$ spanned by all monomials $t_{i_1} \cdots t_{i_k}$ ($i_1, \ldots, i_k \in \{1, \ldots, n\}$), of total degree $k$. The general linear group $\mathrm{GL}_n(K)$ acts naturally on $D_n^1$ and we extend this action so that $\mathrm{GL}_n(K)$ becomes a group of algebra automorphisms of $D_n$. For $f \in D_n \setminus \{0\}$, we consider the following endomorphisms of $R_n$ satisfying the following conditions

$$\begin{aligned}
\alpha_{(f,1)}(y_1) &= y_1 + [[y_1, y_2] \cdot f, [y_3, y_4]], \quad \alpha_{(1,f)}(y_1) = y_1 + [[y_1, y_2], [y_3, y_4] \cdot f], \\
\beta_{(f,1)}(y_1) &= y_1 + [[y_1, y_2] \cdot f, [y_2, y_4]], \quad \gamma_{(f,1)}(y_1) = y_1 + [[y_1, y_2] \cdot f, [y_1, y_4]], \\
\delta_{(f,1)}(y_1) &= y_1 + [[y_1, y_2] \cdot f, [y_2, y_3]], \quad \varepsilon_{(f,1)}(y_1) = y_1 + [[y_2, y_3] \cdot f, [y_2, y_3]] \\
\zeta_{(f,1)}(y_1) &= y_1 + [[y_2, y_3] \cdot f, [y_2, y_4]], \quad \eta_{(f,1)}(y_1) = y_1 + [[y_3, y_4] \cdot f, [y_1, y_3]]
\end{aligned}$$

and $\alpha_{(f,1)}(y_i) = \alpha_{(1,f)}(y_i) = \beta_{(f,1)}(y_i) = \gamma_{(f,1)}(y_i) = \delta_{(f,1)}(y_i) = \zeta_{(f,1)}(y_i) = \eta_{(f,1)}(y_i) = y_i$ for $i \in \{2, \ldots, n\}$. For $n \geq 5$, let $\theta_{(f,1)}$ be the endomorphism of $R_n$ satisfying the conditions $\theta_{(f,1)}(y_1) = y_1 + [[y_2, y_3] \cdot f, [y_4, y_n]]$ and $\theta_{(f,1)}(y_i) = y_i$, $i = 2, \ldots, n$. By Lemma 3.1, the above endomorphisms of $R_n$ are automorphisms of $R_n$. They are determined by the elements $f \in D_n \setminus \{0\}$. By Lemma 3.4, $\mathrm{Ker}\pi_n$ is generated as $K\mathrm{GL}_n(K)$-module by all automorphisms of the form $y_1 \mapsto y_1 + u$, $y_i \mapsto y_i$, $i \geq 2$, $u \in R_n''$. Since every element $u \in R_n''$ is written as a $K$-linear combination of elements of the form $[[w_1, w_2] \cdot f, [w_3, w_4]]$, with $w_1, w_2, w_3, w_4 \in \{y_1, \ldots, y_n\}$ and $f \in D_n \setminus \{0\}$, we obtain $\mathrm{Ker}\pi_n$ is generated as $K\mathrm{GL}_n(K)$-module by all (non-trivial) automorphisms of the form $y_1 \mapsto y_1 + [[w_1, w_2] \cdot f, [w_3, w_4]]$, $y_i \mapsto y_i$, $i \geq 2$, with $w_1, w_2, w_3, w_4 \in \{y_1, \ldots, y_n\}$ and $f \in D_n \setminus \{0\}$. It is easily verified that any element of $\mathrm{Ker}\pi_n$ of the aforementioned form is either of the form $\alpha_{(f,1)}, \alpha_{(1,f)}, \beta_{(f,1)}$ etc. (or their inverses) with $f \in D_n \setminus \{0\}$ or is an automorphism that can be obtained from those by conjugating by suitable tame automorphisms of $R_n$. In particular, for $n, k \geq 4$, $\mathrm{Ker}\overline{\vartheta}_{k,R_n}$ is generated as $K\mathrm{GL}_n(K)$-module by the set $\{\overline{\alpha}_{(1,f)}, \overline{\xi}_{(f,1)} : \overline{\xi}_{(f,1)} \in \{\overline{\alpha}_{(f,1)}, \ldots, \overline{\theta}_{(f,1)}\}; f \in D_n^{k-4} \setminus \{0\}\}$.

**Lemma 3.6** *For positive integers $n, k \geq 4$, $\mathrm{Ker}\overline{\vartheta}_{k,R_n}$ is generated as $K\mathrm{GL}_n(K)$-module by the set $\{\overline{\alpha}_{(f,1)}, \overline{\alpha}_{(1,f)}, \overline{\gamma}_{(f,1)} : f \in D_n^{k-4} \setminus \{0\}\}$.*

*Proof.* Let $\Lambda_{n,k}$ be the $K\mathrm{GL}_n(K)$-submodule of $\mathrm{Ker}\overline{\vartheta}_{k,R_n}$ generated by the set $\{\overline{\alpha}_{(f,1)}, \overline{\alpha}_{(1,f)}, \overline{\gamma}_{(f,1)} : f \in D_n^{k-4} \setminus \{0\}\}$. To prove that $\Lambda_{n,k} = \mathrm{Ker}\overline{\vartheta}_{k,R_n}$, it is enough to show that $\overline{\xi}_{(f,1)} \in$

$\Lambda_{n,k}$ for $\xi_{(f,1)} \in \{\beta_{(f,1)}, \delta_{(f,1)}, \ldots, \theta_{(f,1)}\}$ and $f \in D_n^{k-4} \setminus \{0\}$. Working in Ker$\pi_n$, by direct calculations, we get

$$\beta_{(f,1)} = \alpha^{-1}_{(f,1)}\tau_{3,2}\alpha_{(\tau^{-1}_{3,2}(f),1)}\tau^{-1}_{3,2}, \quad \delta_{(f,1)} = \alpha_{(f,1)}\tau_{4,2}\alpha^{-1}_{(\tau^{-1}_{4,2}(f),1)}\tau^{-1}_{4,2}, \quad \varepsilon_{(f,1)} = \delta_{(f,1)}\tau_{1,3}\delta^{-1}_{(\tau^{-1}_{1,3}(f),1)}\tau^{-1}_{1,3},$$
$$\zeta_{(f,1)} = \beta_{(f,1)}\tau_{1,3}\beta^{-1}_{(\tau^{-1}_{1,3}(f),1)}\tau^{-1}_{1,3}, \quad \eta_{(f,1)} = \alpha_{(1,f)}\tau_{2,3}\alpha^{-1}_{(1,\tau^{-1}_{2,3}(f)}\tau^{-1}_{2,3}$$

and, for $n \geq 5$, $\sigma_{2,4}\theta_{(f,1)}\sigma^{-1}_{2,4} = \alpha_{(1,\sigma_{2,4}(f))}\tau_{1,n}\alpha^{-1}_{(1,\tau^{-1}_{1,n}(\sigma_{2,4}(f)))}\tau^{-1}_{1,n}$. Therefore $\overline{\xi}_{(f,1)} \in \Lambda_{n,k}$ for $\xi_{(f,1)} \in \{\beta_{(f,1)}, \delta_{(f,1)}, \ldots, \theta_{(f,1)}\}$ and $f \in D_n^{k-4} \setminus \{0\}$ and so $\Lambda_{n,k} = \mathrm{Ker}\overline{\vartheta}_{k,R_n}$. $\square$

### 3.3.1 A generating set of $\mathrm{Ker}\overline{\vartheta}_{k,R_n}$ as $K\mathrm{GL}_n(K)$-module

The main aim in this section is to prove the following result.

**Theorem 3.1** *For positive integers $n, k \geq 4$, $\mathrm{Ker}\overline{\vartheta}_{k,R_n}$ is generated as $K\mathrm{GL}_n(K)$-module by the set $\{\overline{\alpha}_{(t_1^{k-4},1)}, \overline{\alpha}_{(1,t_3^{k-4})}\}$.*

For the proof of Theorem 3.1 we need some technical results. For positive integers $n, k \geq 4$, let $\chi_{n,k}$ be the natural $K\mathrm{GL}_n(K)$-isomorphism from $\mathrm{Ker}\overline{\vartheta}_{k,R_n}$ onto $[(R''_n)^k]^{\oplus n}$. Thus, for $\overline{\phi} \in \mathrm{Ker}\overline{\vartheta}_{k,R_n}$, $\chi_{n,k}(\overline{\phi}) = (u_1, \ldots, u_n)$ where $\phi(y_i) = y_i + u_i$, $u_i \in (R''_n)^k$, $i \in \{1, \ldots, n\}$. Throughout this section, we find convenient to work in $[(R''_n)^k]^{\oplus n}$ by means of $\chi_{n,k}$ rather than in $\mathrm{Ker}\overline{\vartheta}_{k,R_n}$. Furthermore, for $i, j, \kappa, \lambda \in \{1, \ldots, n\}$ and $f \in D_n \backslash \{0\}$, we write $u(i, j; f; \kappa, \lambda) = [[y_i, y_j] \cdot f, [y_\kappa, y_\lambda]]$ and $v(i, j; \kappa, \lambda; f) = [[y_i, y_j], [y_\kappa, y_\lambda] \cdot f]$. For any group $G$ and $a, b \in G$, we write $a^b = b^{-1}ab$.

**Proposition 3.1** *For positive integers $n, k \geq 4$, let $N_{n,k,1}$ be the $K\mathrm{GL}_n(K)$-submodule of $\mathrm{Ker}\overline{\vartheta}_{k,R_n}$ generated by $\overline{\alpha}_{(t_1^{k-4},1)}$. Then*

1. *$\overline{\alpha}_{(t_i^{k-4},1)} \in N_{n,k,1}$ for all $i = 2, \ldots, n$.*
2. *$\overline{\alpha}_{(t_1^{m_1}t_2^{m_2},1)} \in N_{n,k,1}$ for all non-negative integers $m_1, m_2$ with $m_1 + m_2 = k - 4$.*
3. *$\overline{\alpha}_{(t_1^{m_1}t_i^{m_i},1)} \in N_{n,k,1}$ for all $i = 2, \ldots, n$ and non-negative integers $m_1, m_i$ with $m_1 + m_i = k - 4$.*
4. *$\overline{\alpha}_{(t_1^{m_1}t_5^{m_5}\cdots t_n^{m_n},1)} \in N_{n,k,1}$ for all non-negative integers $m_1, m_5, \ldots, m_n$ with $m_1 + m_5 + \cdots + m_n = k - 4$.*

5. $\overline{\alpha}_{(t_2^{m_2}t_5^{m_5}\cdots t_n^{m_n},1)} \in N_{n,k,1}$ *for all non-negative integers* $m_2, m_5, \ldots, m_n$ *with* $m_2 + m_5 + \cdots + m_n = k - 4$.

6. $\overline{\alpha}_{(t_1^{m_1}t_2^{m_2}t_5^{m_5}\cdots t_n^{m_n},1)} \in N_{n,k,1}$ *for all non-negative integers* $m_1, m_2, m_5, \ldots, m_n$ *with* $m_1 + m_2 + m_5 + \cdots + m_n = k - 4$.

7. $\overline{\alpha}_{(t_1^{m_1}t_2^{m_2}t_4^{m_4}\cdots t_n^{m_n},1)} \in N_{n,k,1}$ *for all non-negative integers* $m_1, m_2, m_4, \ldots, m_n$ *with* $m_1 + m_2 + m_4 + \cdots + m_n = k - 4$.

8. $\overline{\alpha}_{(t_1^{m_1}t_2^{m_2}t_3^{m_3}\cdots t_n^{m_n},1)} \in N_{n,k,1}$ *for all non-negative integers* $m_1, m_2, m_3, \ldots, m_n$ *with* $m_1 + m_2 + m_3 + \cdots + m_n = k - 4$.

*Proof.* For $k = 4$ our claim is trivially true and so we may assume that $k \geq 5$.

1. Let $i = 2$. It is easily verified that $(\alpha_{(t_1^{k-4},1)})^{\tau_{2,1}^{-1}} = (\alpha_{(t_2^{k-4},1)})^{\sigma_{1,2}}\alpha_{(t_1^{k-4},1)}$. Hence $\overline{\alpha}_{(t_2^{k-4},1)} \in N_{n,k,1}$. Fix $i \geq 3$. For $a \in K\setminus\{0\}$, we have $\tau_{2,i,a}\chi_{n,k}(\overline{\alpha}_{(t_2^{k-4},1)}) \in \chi_{n,k}(N_{n,k,1})$. Now

$$\begin{aligned}\tau_{2,i,a}\chi_{n,k}(\overline{\alpha}_{(t_2^{k-4},1)}) &= ([[y_1, y_2 + ay_i] \cdot (t_2 + at_i)^{k-4}, [y_3, y_4]], 0, \ldots, 0)\\ &= ([[y_1, y_2] \cdot (t_2 + at_i)^{k-4}, [y_3, y_4]], 0, \ldots, 0)+\\ &\quad a([[y_1, y_i] \cdot (t_2 + at_i)^{k-4}, [y_3, y_4]], 0, \ldots, 0)\\ &= \textstyle\sum_{j=0}^{k-4} a^j (\binom{k-4}{j} u(1, 2; t_2^{k-4-j}t_i^j; 3, 4), 0, \ldots, 0)+\\ &\quad \textstyle\sum_{j=1}^{k-3} a^j (\binom{k-4}{j-1} u(1, i; t_2^{k-3-j}t_i^{j-1}; 3, 4), 0, \ldots, 0).\end{aligned}$$

By using a Vandermonde determinant argument (briefly a V.d.a) and taking $j = k - 4$, we have

$$(u(1, 2; t_i^{k-4}; 3, 4), 0, \ldots, 0) + (k - 4)(u(1, i; t_2t_i^{k-5}; 3, 4), 0, \ldots, 0) \in \chi_{n,k}(N_{n,k,1}). \quad (3.3)$$

Since $\tau_{1,i,a}\chi_{n,k}(\overline{\alpha}_{(t_1^{k-4},1)}) \in \chi_{n,k}(N_{n,k,1})$, by similar arguments as above, we have

$$\begin{aligned}\tau_{1,i,a}\chi_{n,k}(\overline{\alpha}_{(t_1^{k-4},1)}) &= \textstyle\sum_{j=0}^{k-4} a^j (\binom{k-4}{j} u(1, 2; t_1^{k-4-j}t_i^j; 3, 4), 0, \ldots, 0)+\\ &\quad \textstyle\sum_{j=1}^{k-3} a^j (\binom{k-4}{j-1} u(i, 2; t_1^{k-3-j}t_i^{j-1}; 3, 4), 0, \ldots, 0).\end{aligned}$$

As before, by using a V.d.a and taking $j = k - 4$, we get

$$(u(1, 2; t_i^{k-4}; 3, 4), 0, \ldots, 0) + (k - 4)(u(i, 2; t_1t_i^{k-5}; 3, 4), 0, \ldots, 0) \in \chi_{n,k}(N_{n,k,1}). \quad (3.4)$$

By using the Jacobi identity in the form $[x,y,z]=[x,z,y]+[x,[y,z]]$, we get

$$u(1,i;t_2t_i^{k-5};3,4)=[[y_1,y_i,y_2]\cdot t_i^{k-5},[y_3,y_4]]=u(1,2;t_i^{k-4};3,4)-u(i,2;t_1t_i^{k-5};3,4)$$

and so (3.3) becomes

$$(k-3)(u(1,2;t_i^{k-4};3,4),0,\ldots,0)-(k-4)(u(i,2;t_1t_i^{k-5};3,4),0,\ldots,0)\in\chi_{n,k}(N_{n,k,1}).$$

Since $k\geq 5$ and $\chi_{n,k}(N_{n,k,1})$ is a vector space, by the above relation and (3.4), we have $(u(1,2;t_i^{k-4};3,4),0,\ldots,0)\in\chi_{n,k}(N_{n,k,1})$. Since $\chi_{n,k}$ is $1-1$, we obtain $\overline{\alpha}_{(t_i^{k-4},1)}\in N_{n,k,1}$ for all $i\in\{3,\ldots,n\}$ and so $\overline{\alpha}_{(t_i^{k-4},1)}\in N_{n,k,1}$ for all $i\in\{2,\ldots,n\}$.

2. For $m_2=0$ or $m_1=0$, our claim follows from the definition of $N_{n,k,1}$ and Proposition 3.1 (1). Thus we may assume that $m_1,m_2>0$. For $a\in K\setminus\{0\}$, we have $\tau_{1,2,a}\chi_{n,k}(\overline{\alpha}_{(t_1^{k-4},1)})\in\chi_{n,k}(N_{n,k,1})$. But

$$\begin{aligned}\tau_{1,2,a}\chi_{n,k}(\overline{\alpha}_{(t_1^{k-4},1)}) &= ([[y_1,y_2]\cdot(t_1+at_2)^{k-4},[y_3,y_4]],0,\ldots,0)\\ &= \textstyle\sum_{j=0}^{k-4}a^j(\binom{k-4}{j}u(1,2;t_1^{k-4-j}t_2^j;3,4),0,\ldots,0).\end{aligned}$$

By using a V.d.a and taking $j=m_2$, we have

$$\binom{k-4}{m_2}(u(1,2;t_1^{m_1}t_2^{m_2};3,4),0,\ldots,0)\in\chi_{n,k}(N_{n,k,1}).$$

Since $\chi_{n,k}(N_{n,k,1})$ is a vector space and $\chi_{n,k}$ is $1-1$, we get $\overline{\alpha}_{(t_1^{m_1}t_2^{m_2},1)}\in N_{n,k,1}$.

3. By the definition of $N_{n,k,1}$ and Proposition 3.1 (1), we may assume that $m_1,m_i>0$ and so $k\geq 6$. For $i=2$, our claim follows from Proposition 3.1 (2). Thus we assume that $i\geq 3$ and fix such an $i$. Since $\overline{\alpha}_{(t_1^{k-4},1)}\in N_{n,k,1}$ and $\tau_{1,i,a}\chi_{n,k}(\overline{\alpha}_{(t_1^{k-4},1)})\in\chi_{n,k}(N_{n,k,1})$ for $a\in K\setminus\{0\}$, we have

$$\textstyle\sum_{j=0}^{k-4}a^j(\binom{k-4}{j}u(1,2;t_1^{k-4-j}t_i^j;3,4),0,\ldots,0)+$$

$$\textstyle\sum_{j=1}^{k-3}a^j(\binom{k-4}{j-1}u(i,2;t_1^{k-3-j}t_i^{j-1};3,4),0,\ldots,0)\in\chi_{n,k}(N_{n,k,1}).$$

By using a V.d.a and taking $j=m_i$, we have

$$\textstyle\binom{k-4}{m_i}(u(1,2;t_1^{m_1}t_i^{m_i};3,4),0,\ldots,0)+$$

$$\textstyle\binom{k-4}{m_i-1}u(i,2;t_1^{m_1+1}t_i^{m_i-1};3,4),0,\ldots,0)\in\chi_{n,k}(N_{n,k,1}).\quad(3.5)$$

By using the Jacobi identity, we have $[y_i, y_2, y_1] = -[y_1, y_i, y_2] + [y_1, y_2, y_i]$ and so

$$u(i,2;t_1^{m_1+1}t_i^{m_i-1};3,4) = u(1,2;t_1^{m_1}t_i^{m_i};3,4) - u(1,i;t_1^{m_1}t_2t_i^{m_i-1};3,4).$$

Therefore (3.5) becomes

$$(\tbinom{k-4}{m_i} + \tbinom{k-4}{m_i-1})(u(1,2;t_1^{m_1}t_i^{m_i};3,4),0,\dots,0)-$$

$$\tbinom{k-4}{m_i-1}u(1,i;t_1^{m_1}t_2t_i^{m_i-1};3,4),0,\dots,0) \in \chi_{n,k}(N_{n,k,1}). \quad (3.6)$$

Since $\overline{\alpha}_{(t_1^{m_1}t_2^{m_i},1)} \in N_{n,k,1}$ and $\tau_{2,i,a}\chi_{n,k}(\overline{\alpha}_{(t_1^{m_1}t_2^{m_i},1)} \in \chi_{n,k}(N_{n,k,1})$ with $a \in K \setminus \{0\}$, we have

$$\textstyle\sum_{j=0}^{m_i} a^j(\binom{m_i}{j}u(1,2;t_1^{m_1}t_2^{m_i-j}t_i^j;3,4),0,\dots,0)+$$
$$\textstyle\sum_{j=1}^{m_i+1} a^j(\binom{m_i}{j-1}u(1,i;t_1^{m_1}t_2^{m_i-j+1}t_i^{j-1};3,4),0,\dots,0) \in \chi_{n,k}(N_{n,k,1}).$$

By using a V.d.a and taking $j = m_i$, we get

$$(u(1,2;t_1^{m_1}t_i^{m_i};3,4),0,\dots,0)+$$

$$\tbinom{m_i}{m_i-1}u(1,i;t_1^{m_1}t_2t_i^{m_i-1};3,4),0,\dots,0) \in \chi_{n,k}(N_{n,k,1}). \quad (3.7)$$

Since $\chi_{n,k}(N_{n,k,1})$ is a vector space and $\chi_{n,k}$ is $1-1$, we obtain, by (3.6) and (3.7), $\overline{\alpha}_{(t_1^{m_1}t_i^{m_i},1)} \in N_{n,k,1}$.

4. We induct on $i$. For $i = 5$, our claim follows from Proposition 3.1 (3). Let $i \geq 5$ and assume that $\overline{\alpha}_{(t_1^{n_1}t_5^{n_5}\cdots t_i^{n_i},1)} \in \chi_{n,k}(N_{n,k,1})$ for all non-negative integers $n_1, n_5, \dots, n_i$ such that $n_1 + n_5 + \cdots + n_i = k-4$. Let $m_1, m_5, \dots, m_{i+1}$ be non-negative integers such that $m_1 + m_5 + \cdots + m_{i+1} = k-4$. Since, by our inductive hypothesis, $\overline{\alpha}_{(t_1^{m_1}t_5^{m_5}\cdots t_i^{m_i+m_{i+1}},1)} \in \chi_{n,k}(N_{n,k,1})$ and $\tau_{i,i+1,a}\chi_{n,k}(\overline{\alpha}_{(t_1^{m_1}t_5^{m_5}\cdots t_i^{m_i+m_{i+1}},1)}) \in \chi_{n,k}(N_{n,k,1})$ with $a \in K \setminus \{0\}$, we have

$$\sum_{j=0}^{m_i+m_{i+1}} a^j(\binom{m_i+m_{i+1}}{j}u(1,2;t_1^{m_1}t_5^{m_5}\dots t_{i-1}^{m_{i-1}}t_i^{m_i+m_{i+1}-j}t_{i+1}^j;3,4),0,\dots,0) \in \chi_{n,k}(N_{n,k,1}).$$

By using a V.d.a and taking $j = m_{i+1}$, we get

$$\binom{m_i+m_{i+1}}{m_{i+1}}(u(1,2;t_1^{m_1}t_5^{m_5}\dots t_i^{m_i}t_{i+1}^{m_{i+1}};0,\dots,0) \in \chi_{n,k}(N_{n,k,1}).$$

By similar arguments given in the proof of Proposition 3.1 (3), we have $\overline{\alpha}_{(t_1^{m_1}t_5^{m_5}\cdots t_{i+1}^{m_{i+1}},1)}$ $\in N_{n,k,1}$ and hence we obtain the result.

5. For $m_2 = 0$, our claim follows from Proposition 3.1 (4) (for $m_1 = 0$). Thus we assume that $m_2 \geq 1$. By Proposition 3.1 (4), $\overline{\alpha}_{(t_1^{m_2} t_5^{m_5} \cdots t_n^{m_n}, 1)} \in N_{n,k,1}$ and so

$$\tau_{1,2,a}\chi_{n,k}(\overline{\alpha}_{(t_1^{m_2} t_5^{m_5} \cdots t_n^{m_n}, 1)}) \in \chi_{n,k}(N_{n,k,1}).$$

But

$$\tau_{1,2,a}\chi_{n,k}(\overline{\alpha}_{(t_1^{m_2} t_5^{m_5} \cdots t_n^{m_n}, 1)}) = \sum_{j=0}^{m_2} a^j (\binom{m_2}{j} u(1, 2; t_1^{m_2-j} t_2^j t_5^{m_5} \cdots t_n^{m_n}; 3, 4), 0, \ldots, 0).$$

By using a V.d.a and taking $j = m_2$, we get $\chi_{n,k}(\overline{\alpha}_{(t_2^{m_2} t_5^{m_5} \cdots t_n^{m_n}, 1)}) \in \chi_{n,k}(N_{n,k,1})$. Since $\chi_{n.k}$ is $1-1$, we have $\overline{\alpha}_{(t_2^{m_2} t_5^{m_5} \cdots t_n^{m_n}, 1)} \in N_{n,k,1}$.

6. For $m_1 = 0$, our claim follows from Proposition 3.1 (5) and so we may assume that $m_1 \geq 1$. By Proposition 3.1 (4), $\overline{\alpha}_{(t_1^{m_1+m_2} t_5^{m_5} \cdots t_n^{m_n}, 1)} \in N_{n,k,1}$. But $\tau_{1,2,a}\chi_{n,k}(\overline{\alpha}_{(t_1^{m_1+m_2} t_5^{m_5} \cdots t_n^{m_n}, 1)})$ $\in \chi_{n,k}(N_{n,k,1})$ with $a \in K \setminus \{0\}$ and

$\tau_{1,2,a}\chi_{n,k}(\overline{\alpha}_{(t_1^{m_1+m_2} t_5^{m_5} \cdots t_n^{m_n}, 1)}) =$

$$\textstyle\sum_{j=0}^{m_1+m_2} a^j (\binom{m_1+m_2}{j} u(1, 2; t_1^{m_1+m_2-j} t_2^j t_5^{m_5} \cdots t_n^{m_n}; 3, 4)., 0, \ldots, 0).$$

By using a V.d.a and taking $j = m_2$, we get

$$\binom{m_1+m_2}{m_2}(u(1, 2; t_1^{m_1} t_2^{m_2} t_5^{m_5} \cdots t_n^{m_n}; 3, 4), 0, \ldots, 0) \in \chi_{n,k}(N_{n,k,1}). \quad (3.8)$$

Since $\chi_{n,k,1}(N_{n,k,1})$ is a vector space and $\chi_{n.k}$ is $1-1$, we have, by (3.8), $\overline{\alpha}_{(t_1^{m_1} t_2^{m_2} t_5^{m_5} \cdots t_n^{m_n}, 1)}$ $\in N_{n,k,1}$.

7. We firstly show that $\overline{\alpha}_{(t_1^{m_1} t_2^{m_2} t_4^{m_4}, 1)} \in N_{n,k,1}$ for all non-negative integers $m_1$, $m_2$, $m_4$ such that $m_1 + m_2 + m_4 = k - 4$. For $m_1 + m_2 = 0$, our claim follows from Proposition 3.1 (1) (for $i = 4$). Thus we may assume that $m_1 + m_2 > 0$. By Proposition 3.1 (3) (for $i = 4$), $\overline{\alpha}_{(t_1^{m_1+m_2} t_4^{m_4}, 1)} \in N_{n,k,1}$. Since $\chi_{n,k}(N_{n,k,1})$ is $K\mathrm{GL}_n(K)$-module, we have $\tau_{1,2,a}\chi_{n,k}(\overline{\alpha}_{(t_1^{m_1+m_2} t_4^{m_4}, 1)}) \in \chi_{n,k}(N_{n,k,1})$ with $a \in K \setminus \{0\}$. But

$$\tau_{1,2,a}\chi_{n,k}(\overline{\alpha}_{(t_1^{m_1+m_2} t_4^{m_4}, 1)}) = \sum_{j=0}^{m_1+m_2} a^j (\binom{m_1+m_2}{j} u(1, 2; t_1^{m_1+m_2-j} t_2^j t_4^{m_4}; 3, 4), 0, \ldots, 0).$$

By using a V.d.a and taking $j = m_2$, we have

$$\binom{m_1+m_2}{m_2}(u(1, 2; t_1^{m_1} t_2^{m_2} t_4^{m_4}; 3, 4), 0, \ldots, 0) \in \chi_{n,k}(N_{n,k,1}). \quad (3.9)$$

Since $\chi_{n,k,1}(N_{n,k,1})$ is a vector space and $\chi_{n.k}$ is $1-1$, we have, by (3.9), $\overline{\alpha}_{(t_1^{m_1}t_2^{m_2}t_4^{m_4},1)} \in N_{n,k,1}$. Next we show that $\overline{\alpha}_{(t_1^{m_1}t_2^{m_2}t_4^{m_4}t_5^{m_5},1)} \in N_{n,k,1}$ for all non-negative integers $m_1$, $m_2$, $m_4$, $m_5$ such that $m_1+m_2+m_4+m_5 = k-4$. Since $\overline{\alpha}_{(t_1^{m_1}t_2^{m_2}t_4^{m_4},1)} \in N_{n,k,1}$, we may assume that $m_5 \geq 1$. Since $\chi_{n,k}(N_{n,k,1})$ is $K\mathrm{GL}_n(K)$-module, we get $\tau_{1,5,a}\chi_{n,k}(\overline{\alpha}_{(t_1^{m_1+m_5}t_2^{m_2}t_4^{m_4},1)}) \in \chi_{n,k}(N_{n,k,1})$ with $a \in K\setminus\{0\}$. For the next few lines, let $w = \tau_{1,5,a}\chi_{n,k}(\overline{\alpha}_{(t_1^{m_1+m_5}t_2^{m_2}t_4^{m_4},1)})$. But

$$\begin{aligned} w &= ([[y_1+ay_5,y_2]\cdot(t_1+at_5)^{m_1+m_5}t_2^{m_2}t_4^{m_4},[y_3,y_4]],0,\ldots,0)\\ &= ([[y_1,y_2]\cdot(t_1+at_5)^{m_1+m_5}t_2^{m_2}t_4^{m_4},[y_3,y_4]],0,\ldots,0)+\\ &\quad a([[y_5,y_2]\cdot(t_1+at_5)^{m_1+m_5}t_2^{m_2}t_4^{m_4},[y_3,y_4]],0,\ldots,0)\\ &= \textstyle\sum_{j=0}^{m_1+m_5} a^j(\binom{m_1+m_5}{j}u(1,2;t_1^{m_1+m_5-j}t_5^jt_2^{m_2}t_4^{m_4};3,4),0,\ldots,0)+\\ &\quad \textstyle\sum_{j=1}^{m_1+m_5+1} a^j(\binom{m_1+m_5}{j-1}u(5,2;t_1^{m_1+m_5+1-j}t_5^{j-1}t_2^{m_2}t_4^{m_4};3,4),0,\ldots,0). \end{aligned}$$

By using a V.d.a and taking $j=m_5$, we have

$$(\tbinom{m_1+m_5}{m_5}u(1,2;t_1^{m_1}t_2^{m_2}t_4^{m_4}t_5^{m_5};3,4),0,\ldots,0)+$$
$$(\tbinom{m_1+m_5}{m_5-1}u(5,2;t_1^{m_1+1}t_2^{m_2}t_4^{m_4}t_5^{m_5-1};3,4),0,\ldots,0) \in \chi_{n,k}(N_{n,k,1}).$$

By using the Jacobi identity in the form $[x,y,z] = [x,z,y]+[x,[y,z]]$, we obtain $[y_5,y_2,y_1] = [y_1,y_2,y_5]-[y_1,y_5,y_2]$ and so

$$((\tbinom{m_1+m_5}{m_5}+\tbinom{m_1+m_5}{m_5-1})(u(1,2;t_1^{m_1}t_2^{m_2}t_4^{m_4}t_5^{m_5};3,4),0,\ldots,0)-$$
$$\tbinom{m_1+m_5}{m_5-1}(u(1,5;t_1^{m_1}t_2^{m_2+1}t_4^{m_4}t_5^{m_5-1};3,4),0,\ldots,0) \in \chi_{n,k}(N_{n,k,1}). \quad (3.10)$$

By Proposition 3.1 (6), $\overline{\alpha}_{(t_1^{m_1}t_2^{m_2+m_5}t_4^{m_4},1)} \in N_{n,k,1}$. Since $\chi_{n,k}(N_{n,k,1})$ is $K\mathrm{GL}_n(K)$-module, we get $\tau_{2,5,a}\chi_{n,k}(\overline{\alpha}_{(t_1^{m_1}t_2^{m_2+m_5}t_4^{m_4},1)}) \in \chi_{n,k}(N_{n,k,1})$. By similar arguments as before, we have

$$\tbinom{m_2+m_5}{m_5}(u(1,2;t_1^{m_1}t_2^{m_2}t_4^{m_4}t_5^{m_5};3,4),0,\ldots,0)+$$
$$\tbinom{m_2+m_5}{m_5-1}u(1,5;t_1^{m_1}t_2^{m_2+1}t_4^{m_4}t_5^{m_5-1};3,4),0,\ldots,0) \in \chi_{n,k}(N_{n,k,1}). \quad (3.11)$$

Since $\chi_{n,k}(N_{n,k,1})$ is a vector space and $\chi_{n,k}$ is $1-1$, we obtain, by (3.10) and (3.11), $\overline{\alpha}_{(t_1^{m_1}t_2^{m_2}t_4^{m_4}t_5^{m_5},1)} \in \chi_{n,k}(N_{n,k,1})$ for all non-negative integers $m_1$, $m_2$, $m_4$, $m_5$ such that $m_1+m_2+m_4+m_5 = k-4$. By a similar inductive argument as in the proof of Proposition 3.1 (4), we get the required result.

8. By Proposition 3.1 (7), $\overline{\alpha}_{(t_1^{m_1}t_2^{m_2}t_4^{m_4+m_3}t_5^{m_5}\cdots t_n^{m_n},1)} \in N_{n,k,1}$ for all non-negative integers $m_1, m_2, \ldots, m_n$ with $m_1+m_2+(m_3+m_4)+m_5+\cdots+m_n = k-4$. For $m_3 = 0$, our claim follows from Proposition 3.1 (7). Thus we assume that $m_3 \geq 1$. Since $\tau_{4,3,a}\chi_{n,k}(\overline{\alpha}_{(t_1^{m_1}t_2^{m_2}t_4^{m_4+m_3}t_5^{m_5}\cdots t_n^{m_n},1)}) \in \chi_{n,k}(N_{n,k,1})$ with $a \in K \setminus \{0\}$, by applying similar arguments as in the proof of Proposition 3.1 (6), we obtain the required result. □

**Corollary 3.1** *For positive integers $n, k \geq 4$, let $N_{n,k,1}$ be the $K\mathrm{GL}_n(K)$-submodule of* $\mathrm{Ker}\overline{\vartheta}_{k,R_n}$ *generated by $\overline{\alpha}_{(t_1^{k-4},1)}$. Then, for all $f \in D_n^{k-4} \setminus \{0\}$, $\overline{\alpha}_{(f,1)} \in N_{n,k,1}$.*

*Proof.* An element $f$ of $D_n^{k-4}\setminus\{0\}$ is (uniquely) written as $K$-linear combination of monomials of the form $t_1^{m_1}\cdots t_n^{m_n}$ with $m_1+\cdots+m_n = k-4$. By Proposition 3.1 (8) and since $N_{n,k,1}$ is a vector space, we obtain the desired result. □

For the proof of the following result, we slightly modify the calculations given in the proof of Proposition 3.1 (1)-(6)–we use $v(i,j;\kappa,\lambda;f)$ instead of $u(i,j;f;\kappa,\lambda)$– and so we avoid the technical details.

**Proposition 3.2** *For positive integers $n, k \geq 4$, let $N_{n,k,2}$ be the $K\mathrm{GL}_n(K)$-submodule of* $\mathrm{Ker}\overline{\vartheta}_{k,R_n}$ *generated by $\overline{\alpha}_{(1,t_3^{k-4})}$. Then*

1. *$\overline{\alpha}_{(1,t_i^{k-4})} \in N_{n,k,2}$ for all $i = 4, \ldots, n$.*

2. *$\overline{\alpha}_{(1,t_3^{m_3}t_4^{m_4})} \in N_{n,k,2}$ for all non-negative integers $m_3, m_4$ with $m_3+m_4 = k-4$.*

3. *$\overline{\alpha}_{(1,t_3^{m_3}t_i^{m_i})} \in N_{n,k,2}$ for all $i = 4, \ldots, n$ and non-negative integers $m_3, m_i$ with $m_3+m_i = k-4$.*

4. *$\overline{\alpha}_{(1,t_3^{m_3}t_5^{m_5}\cdots t_n^{m_n})} \in N_{n,k,2}$ for all non-negative integers $m_3, m_5, \ldots, m_n$ with $m_3+m_5+\cdots+m_n = k-4$.*

5. *$\overline{\alpha}_{(1,t_4^{m_4}t_5^{m_5}\cdots t_n^{m_n})} \in N_{n,k,2}$ for all non-negative integers $m_4, m_5, \ldots, m_n$ with $m_4+m_5+\cdots+m_n = k-4$.*

6. $\overline{\alpha}_{(1,t_3^{m_3}t_4^{m_4}t_5^{m_5}\cdots t_n^{m_n})} \in N_{n,k,2}$ *for all non-negative integers* $m_3, m_4, m_5, \ldots, m_n$ *with* $m_3 + m_4 + m_5 + \cdots + m_n = k - 4$. $\square$

**Proposition 3.3** *For positive integers* $n, k \geq 4$*, let* $N_{n,k} = N_{n,k,1} + N_{n,k,2}$*. Then*

1. $\overline{\gamma}_{(t_1^{m_1}t_2^{m_2}t_5^{m_5}\cdots t_n^{m_n},1)} \in N_{n,k}$ *for all non-negative integers* $m_1, m_2, m_5, \ldots, m_n$ *with* $m_1 + m_2 + m_5 + \cdots + m_n = k - 4$.
2. $\overline{\gamma}_{(t_1^{m_1}t_2^{m_2}t_4^{m_4}\cdots t_n^{m_n},1)} \in N_{n,k}$ *for all non-negative integers* $m_1, m_2, m_4, \ldots, m_n$ *with* $m_1 + m_2 + m_4 + \cdots + m_n = k - 4$.
3. $\overline{\gamma}_{(t_1^{m_1}t_2^{m_2}\cdots t_n^{m_n},1)} \in N_{n,k}$ *for all non-negative integers* $m_1, m_2, \ldots, m_n$ *with* $m_1 + m_2 + \cdots + m_n = k - 4$.

*Proof.*

1. Write $f = t_1^{m_1}t_2^{m_2}t_5^{m_5}\ldots t_n^{m_n}$. By Proposition 3.1 (6), for $a \in K \setminus \{0\}$, we have $\tau_{3,1,a}\chi_{n,k}(\overline{\alpha}_{f,1)}) \in \chi_{n,k}(N_{n,k})$. Hence

   $(u(1,2;f;3,4)+au(1,2;f;1,4), 0, -au(1,2;f;3,4)-a^2u(1,2;f;1,4), 0, \ldots, 0) \in \chi_{n,k}(N_{n,k})$.

   By using a V.d.a, we get

   $$(u(1,2;f;1,4), 0, -u(1,2;f;3,4), 0, \ldots, 0) \in \chi_{n,k}(N_{n,k}). \tag{3.12}$$

   Furthermore, by Proposition 3.2 (6), we have $\sigma_{3,1}\sigma_{2,4}\chi_{n,k}(\overline{\alpha}_{(1,t_3^{m_1}t_4^{m_2}t_5^{m_5}\ldots t_n^{m_n})}) \in \chi_{n,k}(N_{n,k})$ and hence

   $$(0, 0, -u(1,2;f;3,4), 0, \ldots, 0) \in \chi_{n,k}(N_{n,k}). \tag{3.13}$$

   Since $\chi_{n,k}(N_{n,k})$ is a vector space, we have, by (3.12) and (3.13), $\chi_{n,k}(\overline{\gamma}_{(f,1)}) \in \chi_{n,k}(N_{n,k})$. Since $\chi_{n,k}$ is $1-1$, we obtain $\overline{\gamma}_{(f,1)} \in N_{n,k}$.

2. Let $m_1, m_2, m_4, \ldots, m_n$ be non-negative integers such that $m_1 + m_2 + m_4 + \cdots + m_n = k-4$ with $m_4 \geq 1$. By Proposition 3.3 (1), $\overline{\gamma}_{(t_1^{m_1+m_4}t_2^{m_2}t_5^{m_5}\cdots t_n^{m_n},1)} \in N_{n,k}$. For $a \in K\backslash\{0\}$, we have $\tau_{1,4,a}\chi_{n,k}(\overline{\gamma}_{(t_1^{m_1+m_4}t_2^{m_2}t_5^{m_5}\cdots t_n^{m_n},1)}) \in \chi_{n,k}(N_{n,k})$. For the next few lines, we write $\varphi = \tau_{1,4,a}\chi_{n,k}(\overline{\gamma}_{(t_1^{m_1+m_4}t_2^{m_2}t_5^{m_5}\cdots t_n^{m_n},1)})$. But

   $$\begin{aligned}\varphi \;=\; & \textstyle\sum_{j=0}^{m_1+m_4} a^j(\binom{m_1+m_4}{j})u(1,2;t_1^{m_1+m_4-j}t_2^{m_2}t_4^{j}t_5^{m_5}\cdots t_n^{m_n};1,4), 0, \ldots, 0)+ \\ & \textstyle\sum_{j=1}^{m_1+m_4+1} a^j(\binom{m_1+m_4}{j-1})u(4,2;t_1^{m_1+m_4-j+1}t_2^{m_2}t_4^{j-1}t_5^{m_5}\cdots t_n^{m_n};1,4), 0, 0, \ldots, 0).\end{aligned}$$

   By using a V.d.a and taking $j = m_4$, we get

$\binom{m_1+m_4}{m_4}(u(1,2;t_1^{m_1}t_2^{m_2}t_4^{m_4}t_5^{m_5}\cdots t_n^{m_n};1,4),0,\ldots,0)+$

$$\binom{m_1+m_4}{m_4-1}(u(4,2;t_1^{m_1+1}t_2^{m_2}t_4^{m_4-1}t_5^{m_5}\cdots t_n^{m_n};1,4),0,0,\ldots,0)\in\chi_{n,k}(N_{n,k}).$$

By using the Jacobi identity in the form $[x,y,z]=[x,z,y]+[x,[y,z]]$, we have $[y_4,y_2,y_1]=[y_1,y_2,y_4]-[y_1,y_4,y_2]$ and so

$$\left\{\binom{m_1+m_4}{m_4}+\binom{m_1+m_4}{m_4-1}\right\}(u(1,2;t_1^{m_1}t_2^{m_2}t_4^{m_4}t_5^{m_5}\cdots t_n^{m_n};1,4),0,\ldots,0)-$$

$$\binom{m_1+m_4}{m_4-1}(u(1,4;t_1^{m_1}t_2^{m_2+1}t_4^{m_4-1}t_5^{m_5}\cdots t_n^{m_n};1,4),0,\ldots,0)\in\chi_{n,k}(N_{n,k}). \quad (3.14)$$

Furthermore, for $a\in K\setminus\{0\}$, we have $\tau_{2,4,a}\chi_{n,k}(\overline{\gamma}_{(t_1^{m_1}t_2^{m_2+m_4}t_5^{m_5}\cdots t_n^{m_n},1)})\in\chi_{n,k}(N_{n,k})$. We write $\psi=\tau_{2,4,a}\chi_{n,k}(\overline{\gamma}_{(t_1^{m_1}t_2^{m_2+m_4}t_5^{m_5}\cdots t_n^{m_n},1)})$. But

$$\begin{aligned}\psi &= \textstyle\sum_{j=0}^{m_2+m_4}a^j(\binom{m_2+m_4}{j}u(1,2;t_1^{m_1}t_2^{m_2+m_4-j}t_4^{j}t_5^{m_5}\cdots t_n^{m_n};1,4),0,\ldots,0)+\\ &\quad \textstyle\sum_{j=1}^{m_2+m_4+1}a^j(\binom{m_2+m_4}{j-1}u(1,4;t_1^{m_1}t_2^{m_2+m_4-j+1}t_4^{j-1}t_5^{m_5}\cdots t_n^{m_n};1,4),0,0,\ldots,0).\end{aligned}$$

By using a V.d.a and taking $j=m_4$, we get

$\binom{m_2+m_4}{m_4}(u(1,2;t_1^{m_1}t_2^{m_2}t_4^{m_4}t_5^{m_5}\cdots t_n^{m_n};1,4),0,\ldots,0)+$

$$\binom{m_2+m_4}{m_4-1}(u(1,4;t_1^{m_1}t_2^{m_2+1}t_4^{m_4-1}t_5^{m_5}\cdots t_n^{m_n};1,4),0,0,\ldots,0)\in\chi_{n,k}(N_{n,k}). \quad (3.15)$$

Since $\chi_{n,k}(N_{n,k})$ is a vector space, we have, by (3.14) and (3.15), $\chi_{n,k}(\overline{\gamma}_{(t_1^{m_1}t_2^{m_2}t_4^{m_4}\cdots t_n^{m_n},1)})$ $\in\chi_{n,k}(N_{n,k})$. Since $\chi_{n,k}$ is $1-1$, we obtain $\overline{\gamma}_{(t_1^{m_1}t_2^{m_2}t_4^{m_4}\cdots t_n^{m_n},1)}\in N_{n,k}$.

3. Let $m_1,m_2,\ldots,m_n$ be non-negative integers such that $m_1+m_2+\cdots+m_n=k-4$ with $m_3\geq 1$. We assume that $m_3\geq 2$. Similar arguments may be applied for $m_3=1$. By Proposition 3.3 (2), $\overline{\gamma}_{(t_1^{m_1+m_3}t_2^{m_2}t_4^{m_4}\cdots t_n^{m_n},1)}\in N_{n,k}$. For $a\in K\setminus\{0\}$, we have $\tau_{1,3,a}\chi_{n,k}(\overline{\gamma}_{(t_1^{m_1+m_3}t_2^{m_2}t_4^{m_4}\cdots t_n^{m_n},1)})\in\chi_{n,k}(N_{n,k})$. Let $\varphi=\tau_{1,3,a}\chi_{n,k}(\overline{\gamma}_{(t_1^{m_1+m_3}t_2^{m_2}t_4^{m_4}\cdots t_n^{m_n},1)})$. Now

$$\begin{aligned}\varphi &= ([[y_1+ay_3,y_2]\cdot(t_1+at_3)^{m_1+m_3}t_2^{m_2}t_4^{m_4}\cdots t_n^{m_n},[y_1+ay_3,y_4]],0,\ldots,0)\\ &= ([[y_1,y_2]\cdot(t_1+at_3)^{m_1+m_3}t_2^{m_2}t_4^{m_4}\cdots t_n^{m_n},[y_1+ay_3,y_4]],0,\ldots,0)+\\ &\quad a([[y_3,y_2]\cdot(t_1+at_3)^{m_1+m_3}t_2^{m_2}t_4^{m_4}\cdots t_n^{m_n},[y_1+ay_3,y_4]],0,\ldots,0)\\ &= ([[y_1,y_2]\cdot(t_1+at_3)^{m_1+m_3}t_2^{m_2}t_4^{m_4}\cdots t_n^{m_n},[y_1,y_4]],0,\ldots,0)+\\ &\quad a([[y_1,y_2]\cdot(t_1+at_3)^{m_1+m_3}t_2^{m_2}t_4^{m_4}\cdots t_n^{m_n},[y_3,y_4]],0,\ldots,0)+\\ &\quad a([[y_3,y_2]\cdot(t_1+at_3)^{m_1+m_3}t_2^{m_2}t_4^{m_4}\cdots t_n^{m_n},[y_1,y_4]],0,\ldots,0)+\\ &\quad a^2([[y_3,y_2]\cdot(t_1+at_3)^{m_1+m_3}t_2^{m_2}t_4^{m_4}\cdots t_n^{m_n},[y_3,y_4]],0,\ldots,0).\end{aligned}$$

Since $(t_1+at_3)^{m_1+m_3}=\sum_{j=0}^{m_1+m_3} a^j\binom{m_1+m_3}{j}t_1^{m_1+m_3-j}t_3^j$, we have

$$\begin{aligned}\varphi \;=\; & \textstyle\sum_{j=0}^{m_1+m_3} a^j(\binom{m_1+m_3}{j}u(1,2;t_1^{m_1+m_3-j}t_2^{m_2}t_3^jt_4^{m_4}\cdots t_n^{m_n};1,4),0,\ldots,0)+\\ & \textstyle\sum_{j=1}^{m_1+m_3+1} a^j(\binom{m_1+m_3}{j-1}u(1,2;t_1^{m_1+m_3-j+1}t_2^{m_2}t_3^{j-1}t_4^{m_4}\cdots t_n^{m_n};3,4),0,0,\ldots,0)+\\ & \textstyle\sum_{j=1}^{m_1+m_3+1} a^j(\binom{m_1+m_3}{j-1}u(3,2;t_1^{m_1+m_3-j+1}t_2^{m_2}t_3^{j-1}t_4^{m_4}\cdots t_n^{m_n};1,4),0,0,\ldots,0)+\\ & \textstyle\sum_{j=2}^{m_1+m_3+2} a^j(\binom{m_1+m_3}{j-2}u(3,2;t_1^{m_1+m_3-j+2}t_2^{m_2}t_3^{j-2}t_4^{m_4}\cdots t_n^{m_n};3,4),0,0,\ldots,0).\end{aligned}$$

By using a V.d.a and taking $j=m_3$, we get

$$\begin{aligned}&\tbinom{m_1+m_3}{m_3}(u(1,2;t_1^{m_1}t_2^{m_2}t_3^{m_3}t_4^{m_4}\cdots t_n^{m_n};1,4),0,\ldots,0)+\\ &\qquad\tbinom{m_1+m_3}{m_3-1}(u(1,2;t_1^{m_1+1}t_2^{m_2}t_3^{m_3-1}t_4^{m_4}\cdots t_n^{m_n};3,4),0,0,\ldots,0)+\\ &\tbinom{m_1+m_3}{m_3-1}(u(3,2;t_1^{m_1+1}t_2^{m_2}t_3^{m_3-1}t_4^{m_4}\cdots t_n^{m_n};1,4),0,0,\ldots,0)+\\ &\qquad\tbinom{m_1+m_3}{m_3-2}(u(3,2;t_1^{m_1+2}t_2^{m_2}t_3^{m_3-2}t_4^{m_4}\cdots t_n^{m_n};3,4),0,0,\ldots,0)\in\chi_{n,k}(N_{n,k}).\end{aligned}\tag{3.16}$$

By Proposition 3.1 (8),

$$(u(1,2;t_1^{m_1+1}t_2^{m_2}t_3^{m_3-1}t_4^{m_4}\cdots t_n^{m_n};3,4),0,0,\ldots,0)\in\chi_{n,k}(N_{n,k}).\tag{3.17}$$

By Corollary 3.1, $\overline{\alpha}_{(f,1)}\in N_{n,k}$ for all $f\in D_n^{k-4}\backslash\{0\}$. Since $N_{n,k}$ is a $K\mathrm{GL}_n(K)$-module, we have $\overline{\beta}_{(f,1)},\overline{\zeta}_{(f,1)}\in N_{n,k}$ for all $f\in D_n^{k-4}\setminus\{0\}$. Then $\sigma_{2,3}\overline{\zeta}_{(t_1^{m_1+2}t_3^{m_2}t_2^{m_3-2}t_4^{m_4}\cdots t_n^{m_n},1)}\in N_{n,k}$ and so

$$(u(3,2;t_1^{m_1+2}t_2^{m_2}t_3^{m_3-2}t_4^{m_4}\cdots t_n^{m_n};3,4),0,\ldots,0)\in\chi_{n,k}(N_{n,k}).\tag{3.18}$$

By (3.16), (3.17) and (3.18), we have

$$\begin{aligned}&\tbinom{m_1+m_3}{m_3}(u(1,2;t_1^{m_1}t_2^{m_2}t_3^{m_3}t_4^{m_4}\cdots t_n^{m_n};1,4),0,\ldots,0)+\\ &\qquad\tbinom{m_1+m_3}{m_3-1}(u(3,2;t_1^{m_1+1}t_2^{m_2}t_3^{m_3-1}t_4^{m_4}\cdots t_n^{m_n};1,4),0,0,\ldots,0)\in\chi_{n,k}(N_{n,k}).\end{aligned}$$

By using the Jacobi identity in the form $[x,y,z]=[x,z,y]+[x,[y,z]]$, we have $[y_3,y_2,y_1]=[y_1,y_2,y_3]-[y_1,y_3,y_2]$ and so

$$\left\{\tbinom{m_1+m_3}{m_3}+\tbinom{m_1+m_3}{m_3-1}\right\}(u(1,2;t_1^{m_1}t_2^{m_2}t_3^{m_3}t_4^{m_4}\cdots t_n^{m_n};1,4),0,\ldots,0)-$$

$$\binom{m_1+m_3}{m_3-1}(u(1,3;t_1^{m_1}t_2^{m_2+1}t_3^{m_3-1}t_4^{m_4}\cdots t_n^{m_n};1,4),0,\dots,0) \in \chi_{n,k}(N_{n,k}). \quad (3.19)$$

By Proposition 3.3 (2), $\overline{\gamma}_{(t_1^{m_1}t_2^{m_2+m_3}t_4^{m_4}\cdots t_n^{m_n},1)} \in N_{n,k}$. For $a \in K \setminus \{0\}$, we have $\tau_{2,3,a}\chi_{n,k}(\overline{\gamma}_{(t_1^{m_1}t_2^{m_2+m_3}t_4^{m_4}\cdots t_n^{m_n},1)}) \in \chi_{n,k}(N_{n,k})$. Let $\psi = \tau_{2,3,a}\chi_{n,k}(\overline{\gamma}_{(t_1^{m_1}t_2^{m_2+m_3}t_4^{m_4}\cdots t_n^{m_n},1)})$. But

$$\begin{aligned}\psi &= \textstyle\sum_{j=0}^{m_2+m_3} a^j(\binom{m_2+m_3}{j})u(1,2;t_1^{m_1}t_2^{m_2+m_3-j}t_3^{j}t_4^{m_4}\cdots t_n^{m_n};1,4),0,\dots,0)+ \\ &\quad \textstyle\sum_{j=1}^{m_2+m_3+1} a^j(\binom{m_2+m_3}{j-1})u(1,3;t_1^{m_1}t_2^{m_2+m_3-j+1}t_3^{j-1}t_4^{m_4}\cdots t_n^{m_n};1,4),0,0,\dots,0).\end{aligned}$$

By using a V.d.a and taking $j = m_3$, we get

$$\binom{m_2+m_3}{m_3}(u(1,2;t_1^{m_1}t_2^{m_2}t_3^{m_3}t_4^{m_4}\cdots t_n^{m_n};1,4),0,\dots,0)+$$

$$\binom{m_2+m_3}{m_3-1}(u(1,3;t_1^{m_1}t_2^{m_2+1}t_3^{m_3-1}t_4^{m_4}\cdots t_n^{m_n};1,4),0,0,\dots,0) \in \chi_{n,k}(N_{n,k}). \quad (3.20)$$

Since $\chi_{n,k}(N_{n,k})$ is a vector space, we have, by (3.19) and (3.20), $\chi_{n,k}(\overline{\gamma}_{(t_1^{m_1}t_2^{m_2}t_3^{m_3}\cdots t_n^{m_n},1)}) \in \chi_{n,k}(N_{n,k})$. Since $\chi_{n,k}$ is $1-1$, we obtain $\overline{\gamma}_{(t_1^{m_1}t_2^{m_2}t_3^{m_3}\cdots t_n^{m_n},1)} \in N_{n,k}$. □

**Corollary 3.2** *For positive integers $n, k \geq 4$, let $N_{n,k} = N_{n,k,1} + N_{n,k,2}$. Then, for all $f \in D_n^{k-4} \setminus \{0\}$, $\overline{\gamma}_{(f,1)} \in N_{n,k}$.*

*Proof.* An element $f$ of $D_n^{k-4} \setminus \{0\}$ is (uniquely) written as $K$-linear combination of monomials of the form $t_1^{m_1}\cdots t_n^{m_n}$ with $m_1+\cdots+m_n = k-4$. By Proposition 3.3 (3) and since $N_{n,k}$ is a vector space we obtain the desired result. □

**Proposition 3.4** *For positive integers $n, k \geq 4$, let $N_{n,k} = N_{n,k,1} + N_{n,k,2}$. Then, for all $f \in D_n^{k-4} \setminus \{0\}$, $\overline{\alpha}_{(1,f)} \in N_{n,k}$.*

*Proof.* We firstly show that $\overline{\alpha}_{(1,t_1^{m_1}t_3^{m_3}\cdots t_n^{m_n})} \in N_{n,k}$ for all non-negative integers $m_1, m_3, \dots, m_n$ with $m_1 + m_3 + \cdots + m_n = k-4$. By Proposition 3.2 (6), we may assume that $m_1 \geq 1$. Throughout this proof, for $m_1, m_3, \dots, m_n$, we write $\varphi$ for the IA-automorphism of $R_n$ satisfying the conditions $\varphi(y_1) = y_1 + u(2,3;t_1^{m_1}t_2^{m_4}t_3^{m_3}t_5^{m_5}\cdots t_n^{m_n};1,4)$ and $\varphi(y_i) = y_i$, $i \geq 2$. We claim that $\overline{\varphi} \in N_{n,k}$. By Corollary 3.1 and since $N_{n,k}$ is a $K\mathrm{GL}_n(K)$-module, we have $\zeta(f,1) \in N_{n,k}$ for all $f \in D_n^{k-4} \setminus \{0\}$. Notice that, for $j \in \{2, \dots, m_1+2\}$,

$$\sigma_{2,3}\chi_{n,k}(\overline{\zeta}_{(t_1^{m_1+m_3-j+2}t_2^{j-2}t_3^{m_4}t_5^{m_5}\cdots t_n^{m_n},1)}) =$$

$$(u(3,2;t_1^{m_1+m_3-j+2}t_3^{j-2}t_2^{m_4}t_5^{m_5}\cdots t_n^{m_n};3,4),0,\ldots,0)\in\chi_{n,k}(N_{n,k}). \quad (3.21)$$

By Proposition 3.3 (1), $\overline{\gamma}_{(t_1^{m_1+m_3}t_2^{m_4}t_5^{m_5}\cdots t_n^{m_n},1)}\in N_{n,k}$. For $a\in K\setminus\{0\}$, we have

$$\tau_{1,3,a}\chi_{n,k}(\overline{\gamma}_{(t_1^{m_1+m_3}t_2^{m_4}t_5^{m_5}\cdots t_n^{m_n},1)})\in\chi_{n,k}(N_{n,k}).$$

Set $\omega=\tau_{1,3,a}\chi_{n,k}(\overline{\gamma}_{(t_1^{m_1+m_3}t_2^{m_4}t_5^{m_5}\cdots t_n^{m_n},1)})$. By using similar arguments as in the proof of Proposition 3.3 (3), we have

$$\begin{aligned}\omega \; = \; & \textstyle\sum_{j=0}^{m_1+m_3} a^j(\binom{m_1+m_3}{j}u(1,2;t_1^{m_1+m_3-j}t_3^{j}t_2^{m_4}t_5^{m_5}\cdots t_n^{m_n};1,4),0,\ldots,0)+\\ & \textstyle\sum_{j=1}^{m_1+m_3+1} a^j(\binom{m_1+m_3}{j-1}u(1,2;t_1^{m_1+m_3-j+1}t_3^{j-1}t_2^{m_4}t_5^{m_5}\cdots t_n^{m_n};3,4),0,\ldots,0)+\\ & \textstyle\sum_{j=1}^{m_1+m_3+1} a^j(\binom{m_1+m_3}{j-1}u(3,2;t_1^{m_1+m_3-j+1}t_3^{j-1}t_2^{m_4}t_5^{m_5}\cdots t_n^{m_n};1,4),0,\ldots,0)+\\ & \textstyle\sum_{j=2}^{m_1+m_3+2} a^j(\binom{m_1+m_3}{j-2}u(3,2;t_1^{m_1+m_3-j+2}t_3^{j-2}t_2^{m_4}t_5^{m_5}\cdots t_n^{m_n};3,4),0,\ldots,0).\end{aligned}$$

Having in mind (3.21) and the definitions of $\overline{\alpha}_{(f,1)}$ and $\overline{\gamma}_{(f,1)}$, the above expression becomes

$$\begin{aligned}\omega \; = \; & \textstyle\sum_{j=0}^{m_1+m_3} a^j\binom{m_1+m_3}{j}\chi_{n,k}(\overline{\gamma}_{(t_1^{m_1+m_3-j}t_3^{j}t_2^{m_4}t_5^{m_5}\cdots t_n^{m_n};1)})+\\ & \textstyle\sum_{j=1}^{m_1+m_3+1} a^j\binom{m_1+m_3}{j-1}\chi_{n,k}(\overline{\alpha}_{(t_1^{m_1+m_3-j+1}t_3^{j-1}t_2^{m_4}t_5^{m_5}\cdots t_n^{m_n},1)})+\\ & \textstyle\sum_{j=1}^{m_1+m_3+1} a^j(\binom{m_1+m_3}{j-1}u(3,2;t_1^{m_1+m_3-j+1}t_3^{j-1}t_2^{m_4}t_5^{m_5}\cdots t_n^{m_n};1,4),0,\ldots,0)+\\ & \textstyle\sum_{j=2}^{m_1+m_3+2} a^j\binom{m_1+m_3}{j-2}(\sigma_{2,3}\chi_{n,k}(\overline{\zeta}_{(t_1^{m_1+m_3-j+2}t_2^{j-2}t_3^{m_4}t_5^{m_5}\cdots t_n^{m_n},1)})).\end{aligned}$$

By Corollary 3.1, Corollary 3.2 and since $\chi_{n,k}(N_{n,k})$ is a $K\mathrm{GL}_n(K)$-module, we have

$$\sum_{j=1}^{m_1+m_3+1} a^j(\binom{m_1+m_3}{j-1}u(3,2;t_1^{m_1+m_3-j+1}t_3^{j-1}t_2^{m_4}t_5^{m_5}\cdots t_n^{m_n};1,4),0,\ldots,0)\in\chi_{n,k}(N_{n,k}).$$

By using a V.d.a and taking $j=m_3+1$, we obtain

$$\binom{m_1+m_3}{m_3}(u(3,2;t_1^{m_1}t_3^{m_3}t_2^{m_4}t_5^{m_5}\cdots t_n^{m_n};1,4),0,\ldots,0)\in\chi_{n,k}(N_{n,k}).$$

Since $\chi_{n,k}(N_{n,k})$ is a vector space and $\chi_{n,k}$ is $1-1$, we have $\overline{\varphi}\in N_{n,k}$. Since $\sigma_{2,4}\chi_{n,k}(\overline{\varphi})=\chi_{n,k}(\overline{\alpha}_{(1,t_1^{m_1}t_3^{m_3}t_4^{m_4}\cdots t_n^{m_n})})$ and $\chi_{n,k}$ is $1-1$, we have $\overline{\alpha}_{(1,t_1^{m_1}t_3^{m_3}t_4^{m_4}\cdots t_n^{m_n})}\in N_{n,k}$ for all non-negative integers $m_1,m_3,\ldots,m_n$ with $m_1+m_3+\cdots+m_n=k-4$.

Next, we show that $\overline{\alpha}_{(1,f)}\in N_{n,k}$ for all $f\in D_n^{k-4}\setminus\{0\}$. For $a\in K\setminus\{0\}$, we have $\tau_{1,2,a}\chi_{n,k}(\overline{\alpha}_{(1,t_1^{m_1+m_2}t_3^{m_3}t_4^{m_4}\cdots t_n^{m_n})})\in\chi_{n,k}(N_{n,k})$. Clearly we may assume that $m_2\geq 1$. Now

$$\tau_{1,2,a}\chi_{n,k}(\overline{\alpha}_{(1,t_1^{m_1+m_2}t_3^{m_3}t_4^{m_4}\cdots t_n^{m_n})})=$$

$$\textstyle\sum_{j=0}^{m_1+m_2} a^j(\binom{m_1+m_2}{j}(u(1,2;3,4;t_1^{m_1+m_2-j}t_2^{j}t_3^{m_3}\cdots t_n^{m_n}),0,\ldots,0)).$$

By using a V.d.a and taking $j = m_2$, we get

$$\binom{m_1+m_2}{m_2}(u(1,2;3,4;t_1^{m_1}t_2^{m_2}t_3^{m_3}\cdots t_n^{m_n}),0,\ldots,0) \in \chi_{n,k}(N_{n,k}).$$

Since $N_{n,k}$ is a vector space and $\chi_{n,k}$ is $1-1$, we get

$$\overline{\alpha}_{(1,t_1^{m_1}t_2^{m_2}t_3^{m_3}t_4^{m_4}\cdots t_n^{m_n})} \in N_{n,k} \tag{3.22}$$

for all non-negative integers $m_1,\ldots,m_n$ with $m_1+\cdots+m_n = k-4$. An element $f$ of $D_n^{k-4}\setminus\{0\}$ is (uniquely) written as $K$-linear combination of monomials of the form $t_1^{m_1}\cdots t_n^{m_n}$ with $m_1+\cdots+m_n = k-4$. By (3.22) and since $N_{n,k}$ is a vector space we obtain, for all $f \in D_n^{k-4}\setminus\{0\}$, $\overline{\alpha}_{(1,f)} \in N_{n,k}$. □

*Proof of Theorem 3.1.* Let $n,k$ be positive integers, with $n,k \geq 4$. By Lemma 3.6, $\text{Ker}\overline{\vartheta}_{k,R_n}$ is generated as $K\text{GL}_n(K)$-module by the set $\{\overline{\alpha}_{(f,1)},\overline{\alpha}_{(1,f)},\overline{\gamma}_{(f,1)} : f \in D_n^{k-4}\setminus\{0\}\}$. To prove our claim, it is enough to show that $\text{Ker}\overline{\vartheta}_{k,R_n} \subseteq N_{n,k}$. By Corollary 3.1, Corollary 3.2 and Proposition 3.4, we obtain the desired result. □

*Proof of Theorem 1.1.* For a positive integer $\kappa$, we write $\omega_\kappa = \alpha_{(t_1^{\kappa-1},1)}$. Let $G$ be the subgroup of $\text{Aut}(R_n)$ generated by the tame automorphisms $T_{R_n}$ and $\omega_1,\omega_2,\ldots$. To prove that $G$ is dense in $\text{Aut}(R_n)$ with respect to the formal power series topology, it is enough to show that $\bar{\text{I}}_kG = \bar{\text{I}}_k\text{A}_{R_n}$ for all $k \geq 2$. The natural mapping from $R_n$ onto $M_n$ induces a group homomorphism $\pi_n$ from $\text{Aut}(R_n)$ into $\text{Aut}(M_n)$. By Lemma 3.1, $\pi_n$ is surjective. Since $R_n/R_n' \cong M_n/M_n' \cong L_n/L_n'$ as $K$-vector spaces, the restriction $\widetilde{\pi}_n$ of $\pi_n$ on $\text{IA}(R_n)$ is a group epimorphism from $\text{IA}(R_n)$ onto $\text{IA}(M_n)$. Furthermore, for all $k \geq 2$, $\widetilde{\pi}_n$ induces a group epimorphism $\widetilde{\pi}_{n,k}$ from $\text{I}_k\text{A}_{R_n}$ onto $\text{I}_k\text{A}_{M_n}$. Clearly $\widetilde{\pi}_{n,k}$ induces a vector space epimorphism, denoted $\widetilde{\pi}_{n,k}$ as well, from $\bar{\text{I}}_kA_{R_n}$ onto $\bar{\text{I}}_kA_{M_n}$. Recall that $\overline{\vartheta}_{k,R_n}$ is a $K\text{GL}_n(K)$-module epimorphism from $\bar{\text{I}}_k\text{E}_{R_n}$ onto $\bar{\text{I}}_k\text{E}_{M_n}$. In particular if $\overline{\phi}_R \in \bar{\text{I}}_k\text{E}_{R_n}$, with $\phi_R \in \text{I}_k\text{E}_{R_n}$ satisfying the conditions $\phi_R(y_i) = y_i + u_i$, $u_i \in R_n^k$, $i = 1,\ldots,n$, then $\overline{\vartheta}_{k,R_n}(\overline{\phi}_R) = \overline{\phi}_M$, where $\phi_M \in \text{I}_k\text{E}_{M_n}$ such that $\phi_M(z_i) = z_i + \mu_{k,R_n}(u_i)$, $i = 1,\ldots,n$. Moreover, for all $g \in \text{GL}_n(K)$, $\overline{\vartheta}_{k,R_n}(\overline{g\phi_R g^{-1}}) = \overline{g\phi_M g^{-1}}$. By Lemma 3.2, $\text{Ker}\overline{\vartheta}_{k,R_n} \subseteq \bar{\text{I}}_k\text{A}_{R_n}$. It is clear enough that $\text{Ker}\widetilde{\pi}_{n,k} = \text{Ker}\overline{\vartheta}_{k,R_n}$. We claim that $\text{Ker}\overline{\vartheta}_{k,R_n} \subseteq \bar{\text{I}}_kG$ for all $k \geq 4$. By Theorem 3.1, it is enough to show that $\overline{\alpha}_{(1,t_3^{k-4})} \in \bar{\text{I}}_kG$ for all $k \geq 4$. For $k = 4$, our claim is trivially true. Thus we may assume that $k \geq 5$. Let $\phi,\phi_1$ be the tame automorphisms of $R_n$ satisfying the conditions $\phi(y_2) = y_2 + [y_3,y_4]$, $\phi_1(y_1) = y_1 + [y_2,y_3]$, $\phi(y_i) = y_i$ for $i \neq 2$ and $\phi_1(y_j) = y_j$ for

$j \geq 2$. By direct calculations, we get

$$\beta_{(1,1)}^{-1}\phi\beta_{(1,1)}\phi^{-1} = (\sigma_{3,4}\alpha_{(1,t_3)}^{-1}\sigma_{3,4}^{-1})(\sigma_{2,4}\beta_{(t_3,1)}\sigma_{2,4}^{-1})(\sigma_{2,3}\alpha_{(t_4,1)}\sigma_{2,3}^{-1}) \tag{3.23}$$

and, for $k \geq 6$,

$$\gamma_{(t_2^{k-5},1)}^{-1}\phi_1^{-1}\gamma_{(t_2^{k-5},1)}\phi_1 = \sigma_{2,4}\sigma_{3,4}\alpha_{(1,t_3^{k-4})}\sigma_{3,4}^{-1}\sigma_{2,4}^{-1}. \tag{3.24}$$

Notice that $N_{n,k,1} \subseteq \bar{\mathrm{I}}_k G$. By Corollary 3.1, $\overline{\alpha}_{(f,1)} \in \bar{\mathrm{I}}_k G$ for all $f \in D_n^{k-4} \setminus \{0\}$. Since $\beta_{(f,1)} = \alpha_{(f,1)}^{-1}\tau_{2,3}\alpha_{(\tau_{3,2}^{-1}(f),1)}\tau_{3,2}^{-1}$ and $G$ is a $T_{R_n}$-invariant group, we obtain $\overline{\beta}_{(f,1)} \in \bar{\mathrm{I}}_k G$ for all $f \in D_n^{k-4} \setminus \{0\}$. By (3.23), we have $\overline{\alpha}_{(1,t_3)} \in \bar{\mathrm{I}}_5 G$ and so $\mathrm{Ker}\overline{\vartheta}_{5,R_n} \subseteq \bar{\mathrm{I}}_5 G$. Since $N_{n,5} \subseteq \bar{\mathrm{I}}_5 G$, we have, by Corollary 3.2, $\overline{\gamma}_{(f,1)} \in \bar{\mathrm{I}}_5 G$ for all $f \in D_n^1 \setminus \{0\}$. Since $\overline{[\gamma_{(f,1)}, \phi_1]} \in \bar{\mathrm{I}}_6 G$ and $\bar{\mathrm{I}}_6 G$ is a $K\mathrm{GL}_n(K)$-module, we get, by (3.24) (for $k = 6$), $\overline{\alpha}_{(1,t_3^2)} \in \bar{\mathrm{I}}_6 G$ and so $\mathrm{Ker}\overline{\vartheta}_{6,R_n} \subseteq \bar{\mathrm{I}}_6 G$. Proceeding in this way, we obtain $\mathrm{Ker}\overline{\vartheta}_{k,R_n} \subseteq \bar{\mathrm{I}}_k G$ for all $k \geq 5$. Therefore $\mathrm{Ker}\overline{\vartheta}_{k,R_n} \subseteq \bar{\mathrm{I}}_k G$ for all $k \geq 4$.

Let $G^* = \pi_n(G)$. Since $G^* = T_{M_n}$, we have, by a result of Bryant and Drensky [5, Theorem 3.11 (i)], $G^*$ is dense in $\mathrm{Aut}(M_n)$. Therefore, for all $k \geq 2$, $\bar{\mathrm{I}}_k G^* = \bar{\mathrm{I}}_k \mathrm{A}_{M_n}$. Since, for $k = 2, 3$, $\bar{\mathrm{I}}_k \mathrm{A}_{R_n} \cong \bar{\mathrm{I}}_k \mathrm{A}_{M_n} = \bar{\mathrm{I}}_k T_{M_n}$ as $K\mathrm{GL}_n(K)$-modules, we may assume that $k \geq 4$. Since $\widetilde{\pi}_{n,k}(\bar{\mathrm{I}}_k \mathrm{A}_{R_n}) = \bar{\mathrm{I}}_k \mathrm{A}_{M_n}$ and $\mathrm{Ker}\widetilde{\pi}_{n,k} = \mathrm{Ker}\overline{\vartheta}_{k,R_n}$, we have $\bar{\mathrm{I}}_k \mathrm{A}_{R_n}/\mathrm{Ker}\overline{\vartheta}_{k,R_n} \cong \bar{\mathrm{I}}_k \mathrm{A}_{M_n}$. By definitions, $\bar{\mathrm{I}}_k G \subseteq \bar{\mathrm{I}}_k \mathrm{A}_{R_n}$ for all $k \geq 2$. Since $\widetilde{\pi}_{n,k}(\bar{\mathrm{I}}_k G) = \bar{\mathrm{I}}_k G^*$ and $\mathrm{Ker}\overline{\vartheta}_{k,R_n} \subseteq \bar{\mathrm{I}}_k G$, we get $\bar{\mathrm{I}}_k G/\mathrm{Ker}\overline{\vartheta}_{k,R_n} \cong \bar{\mathrm{I}}_k G^*$. Since $\bar{\mathrm{I}}_k \mathrm{A}_{M_n} = \bar{\mathrm{I}}_k G^*$, we have $\bar{\mathrm{I}}_k \mathrm{A}_{R_n}/\mathrm{Ker}\overline{\vartheta}_{k,R_n} \cong \bar{\mathrm{I}}_k G/\mathrm{Ker}\overline{\vartheta}_{k,R_n}$ as vector spaces for all $k \geq 4$. Since $\dim(\bar{\mathrm{I}}_k \mathrm{A}_{R_n})$ is finite, $\dim(\bar{\mathrm{I}}_k G) = \dim(\bar{\mathrm{I}}_k \mathrm{A}_{R_n})$ and $\bar{\mathrm{I}}_k G \subseteq \bar{\mathrm{I}}_k \mathrm{A}_{R_n}$, we obtain $\bar{\mathrm{I}}_k G = \bar{\mathrm{I}}_k \mathrm{A}_{R_n}$ for all $k \geq 4$. Since $\bar{\mathrm{I}}_k G = \bar{\mathrm{I}}_k \mathrm{A}_{R_n}$ for $k = 2, 3$, we have $\bar{\mathrm{I}}_k G = \bar{\mathrm{I}}_k \mathrm{A}_{R_n}$ for all $k \geq 2$ and so $G$ is dense in $\mathrm{Aut}(R_n)$ with respect to the formal power series topology. □

# 4 Non-tame automorphisms

For a positive integer $n$, with $n \geq 4$, let $A_n = K\langle x_1, \ldots, x_n\rangle$ be the free associative algebra over $K$ (with identity element) freely generated by the variables $x_1, \ldots, x_n$. For a non-negative integer $m$, let $A_n^m$ be the subspace of $A_n$ spanned by all monomials of degree $m$ with the convention that 1 is the only monomial of degree 0. Thus $A_n = \bigoplus_{m \geq 0} A_n^m$ with $A_n^0 = K$. Give $A_n$ the structure of a Lie algebra by defining the Lie bracket $[u, v] = uv - vu$ for all $u, v \in A_n$. The Lie subalgebra of $A_n$ generated by the set $\mathcal{A} = \{x_1, \ldots, x_n\}$ is a free Lie algebra on $\mathcal{A}$ (see [2], [4], [18]). Thus we consider the free Lie algebra $L_n$ embedded into $A_n$. Throughout this section, we write $C_n = L_n/J_n$ where $J_n = \gamma_3(L'_n)$. For $j \in \{1, \ldots, n\}$, let $q_j = x_j + J_n$ and so the set $\{q_1, \ldots, q_n\}$ is a free generating set of $C_n$. For a positive integer $\kappa$, let $\varphi_\kappa$ be the IA-automorphism of $C_n$ satisfying the conditions $\varphi_\kappa(q_1) = q_1 + [q_1, q_2, \ {}_{(\kappa-1)}q_1, [q_3, q_4]]$ and $\varphi_\kappa(q_j) = q_j$, $j \geq 2$. We show that each $\varphi_\kappa$ is a non-tame automorphism of $C_n$.

Every element $f \in A_n$ can be written uniquely in the form $f = f_0 + \sum x_i f_i$, where $f_0 \in K$ and $f_i \in A_n^i$ with $i \geq 1$. Following [17] (see, also [5]) we define partial derivatives $\partial_i$ by $\partial_i(f) = f_i$. For $f, g \in A_n$ and $\alpha, \beta \in K$,

$$\begin{aligned} \partial_i(\alpha f + \beta g) &= \alpha\partial_i(f) + \beta\partial_i(g), \\ \partial_i(fg) &= \partial_i(f)g + \varepsilon(f)\partial_i(g) \qquad \text{and} \\ \partial_i([f,g]) &= \partial_i(f)g + \varepsilon(f)\partial_i(g) - \partial_i(g)f - \varepsilon(g)\partial_i(f) \end{aligned}$$

where $\varepsilon : A_n \to K$ is the augmentation homomorphism with kernel $\Delta$. The Jacobian matrix of an endomorphism $\phi$ of $A_n$ is $J(\phi) = (\partial_i(\phi(x_j)))$. Since $L_n \subset A_n$, we can define $J(\phi)$ also for $\phi \in \mathrm{End}(L_n)$. An element $f \in A_n$ is balanced [6, Section 3] (see, also [8, Section 2]) if $f$ belongs to the vector subspace of $A_n^m$ spanned by all elements of the form $x_{i_1}x_{i_2}\cdots x_{i_m} - x_{i_2}\cdots x_{i_m}x_{i_1}$. Let $\phi$ be an endomorphism of $L_n$ such that $\phi(x_i) \equiv x_i + f_i$ modulo $\gamma_m(L_n)$, $i \in \{1, \ldots, n\}$, where for each $i$, $f_i$ is homogeneous of degree $m$ or $f_i = 0$. As shown in [6, Theorem 3.7], if $\phi$ is an automorphism of $L_n$, then $\partial_1(f_1) + \cdots + \partial_n(f_n)$ is balanced.

**Lemma 4.1** *Let $m$ be a positive integer and $j \in \{2, \ldots, n\}$. Then*

1. $[x_1, x_j, \ {}_m x_1] = \left(\sum_{k=1}^{m+1}(-1)^{k-1}\binom{m+1}{k}x_1^k x_j x_1^{m+1-k}\right) - x_j x_1^{m+1}$.
2. $\partial_j([x_1, x_j, \ {}_m x_1]) = -x_1^{m+1}$.
3. $\partial_1([x_1, x_2, \ {}_m x_1, [x_3, x_4]]) = \left(\sum_{k=1}^{m+1}(-1)^{k-1}\binom{m+1}{k}x_1^{k-1}x_2 x_1^{m+1-k}\right)[x_3, x_4]$.

*Proof.* We induct on $m$. Let $m = 1$. Then

$$\begin{aligned} [x_1, x_j, x_1] &= [x_1, x_j]x_1 - x_1[x_1, x_j] \\ &= (2x_1x_jx_1 - x_1^2x_j) - x_jx_1^2 \\ &= \left(\sum_{k=1}^{2}(-1)^{k-1}\binom{2}{k}x_1^k x_j x_1^{2-k}\right) - x_j x_1^2. \end{aligned}$$

By our inductive hypothesis,

$$[x_1, x_j, \ {}_m x_1] = \left(\sum_{k=1}^{m+1}(-1)^{k-1}\binom{m+1}{k}x_1^k x_j x_1^{m+1-k}\right) - x_j x_1^{m+1}.$$

Thus

$$\begin{aligned}[x_1, x_j, {}_{(m+1)}x_1] &= [x_1, x_j, {}_m x_1]x_1 - x_1[x_1, x_j, {}_m x_1] \\ &= (m+2)x_1 x_j x_1^{m+1} + \textstyle\sum_{k=2}^{m+1}(-1)^{k-1}(\binom{m+1}{k} + \binom{m+1}{k-1})x_1^k x_j x^{m+2-k} + \\ &\quad (-1)^{m+1}x_1^{m+2}x_j) - x_j x_1^{m+2} \\ &= \left(\textstyle\sum_{k=1}^{m+2}(-1)^{k-1}\binom{m+2}{k}x_1^k x_j x_1^{m+2-k}\right) - x_j x_1^{m+2}.\end{aligned}$$

Therefore we obtain the required result. By the definition of the partial derivative and Lemma 4.1 (1), we have $\partial_j([x_1, x_j, {}_m x_1]) = -x_1^{m+1}$. Let $v = [x_1, x_2, {}_m x_1, [x_3, x_4]]$. Since $\partial_1(v) = \partial_1([x_1, x_2, {}_m x_1])[x_3, x_4]$, by the definition of the partial derivative and Lemma 4.1 (1) (for $j = 2$), we have the required result. □

Let $M_2(K[z])$ be the associative $K$-algebra of $2 \times 2$ matrices over $K[z]$. Let $\Phi$ be the mapping from $\mathcal{A}$ into $M_2(K[z])$ such that $\Phi(x_i) = \mathrm{Id}_2 = \begin{pmatrix} 1 & 0 \\ 0 & 1 \end{pmatrix}$, with $i \in \{1, \ldots, n\} \setminus \{2, 3, 4\}$, $\Phi(x_2) = \begin{pmatrix} z & -z \\ z & -z \end{pmatrix}$, $\Phi(x_3) = \begin{pmatrix} 0 & z \\ 0 & 0 \end{pmatrix}$ and $\Phi(x_4) = \begin{pmatrix} 0 & 0 \\ z & 0 \end{pmatrix}$. Since $A_n$ is a free associative $K$-algebra on $\mathcal{A}$, $\Phi$ extends to an algebra homomorphism from $A_n$ into $M_2(K[z])$. For a positive integer $m$, let $\Lambda^{(m)}$ be the ideal of $K[z]$ generated by $z^m$. Observe that if $i_1, \ldots, i_\beta, j_1, \ldots, j_\gamma \in \{1, \ldots, n\}$ with $\beta, \gamma \geq 2$, then the entries of $\Phi([x_{i_1}, \ldots, x_{i_\beta}])\Phi([x_{j_1}, \ldots, x_{j_\gamma}])$ belong to $\Lambda^{(4)}$. Fix $i \in \{1, \ldots, n\}$ and assume that $\partial_i([x_{k_1}, \ldots, x_{k_\delta}]) \neq 0$ where $k_1, \ldots, k_\delta \in \{1, \ldots, n\}$ and $\delta \geq 2$. By using the Lie bracket $[u, v] = uv - vu$, each Lie commutator $[x_{k_1}, \ldots, x_{k_\delta}]$ is a $\mathbb{Z}$-linear combination of monomials $x_{j_1} \cdots x_{j_\delta}$ with $j_1, \ldots, j_\delta$ are $k_1, \ldots, k_\delta$ in some order. Write $[x_{k_1}, \ldots, x_{k_\delta}] = \sum_{j=1}^{n} x_j f_j(x_{k_1}, \ldots, x_{k_\delta})$ where $f_j(x_{k_1}, \ldots, x_{k_\delta}) \in A_n^{\delta-1}$. By the definition of the partial derivative, $\partial_i([x_{k_1}, \ldots, x_{k_\delta}]) = f_i(x_{k_1}, \ldots, x_{k_\delta})$. So $\Phi(\partial_i([x_{k_1}, \ldots, x_{k_\delta}])) = f_i(\Phi(x_{k_1}), \ldots, \Phi(x_{k_\delta}))$. Notice that if each $\Phi(x_{k_e}) = \mathrm{Id}_2$ ($e \in \{1, \ldots, \delta\}$), then $\Phi(\partial_i([x_{k_1}, \ldots, x_{k_\delta}])) = \alpha\ \mathrm{Id}_2$ with $\alpha \in K$. Otherwise the entries of $\Phi(\partial_i([x_{k_1}, \ldots, x_{k_\delta}]))$ belong to $\Lambda^{(1)}$.

**Proposition 4.1** *For a positive integer $\kappa$, let $\varphi_\kappa$ be the automorphism of $C_n$ satisfying the conditions $\varphi_\kappa(q_1) = q_1 + [q_1, q_2, {}_{(\kappa-1)}q_1, [q_3, q_4]]$ and $\varphi_\kappa(q_i) = q_i$, $i = 2, \ldots, n$. Then each $\varphi_\kappa$ is non-tame.*

*Proof.* Let $\kappa = 1$. Since $J_n \subseteq \gamma_5(L_n)$, the automorphism $\varphi_1$ is non-tame (see [10, Lemma 2.1]). Thus we may assume that $\kappa \geq 2$. Throughout this proof, let $v_\kappa = [x_1, x_2, {}_{(\kappa-1)}x_1, [x_3, x_4]]$.

Suppose that $\varphi_\kappa$ is induced by an automorphism $\psi_\kappa$ of $L_n$. Then we can write $\psi_\kappa(x_1) = x_1 + v_\kappa + u_1$ and $\psi_\kappa(x_i) = x_i + u_i$ where each $u_i$ belongs to $J_n$. Let $\kappa = 2$. Since $J_n \subseteq \gamma_6(L_n)$, we have $\psi_2(x_1) \equiv x_1 + v_2$ modulo $\gamma_6(L_n)$ and $\psi_2(x_i) \equiv x_i$ modulo $\gamma_6(L_n)$. By Lemma 4.1 (3) (for $m = 1$), we get $\partial_1(v_2) = 2x_2x_1[x_3, x_4] - x_1x_2[x_3, x_4]$ being not balanced and so $\psi_2$ is not an automorphism. From now on we assume that $\kappa \geq 3$. By Lemma 4.1 (3) (for $m = \kappa - 1$), we have

$$\partial_1(v_\kappa) = \left(\sum_{\lambda=1}^{\kappa}(-1)^{\lambda-1}\binom{\kappa}{\lambda}x_1^{\lambda-1}x_2x_1^{\kappa-\lambda}\right)[x_3, x_4]. \tag{4.1}$$

For $B \in M_2(K[z])$, we write $\mathrm{Tr}(B)$ for the trace of $B$. The map $\mathrm{Tr} : M_2(K[z]) \to K[z]$ is a $K$-linear mapping and $\mathrm{Tr}(B_1B_2 - B_2B_1) = 0$. If $u = x_{i_1}x_{i_2}\cdots x_{i_m} - x_{i_2}\cdots x_{i_m}x_{i_1}$, with $m \geq 2$, then

$$\Phi(u) = \Phi(x_{i_1})\Phi(x_{i_2})\cdots\Phi(x_{i_m}) - \Phi(x_{i_2})\cdots\Phi(x_{i_m})\Phi(x_{i_1}) = B_{i_1}(B_{i_2}\cdots B_{i_m}) - (B_{i_2}\cdots B_{i_m})B_{i_1},$$

with $\Phi(x_{i_j}) = B_{i_j}$, $j \in \{1, \ldots, n\}$ and so $\mathrm{Tr}(\Phi(u)) = 0$. Let $a_1 = [x_{i_1}, \ldots, x_{i_{m(1)}}]$, $a_2 = [x_{j_1}, \ldots, x_{j_{m(2)}}]$ and $a_3 = [x_{k_1}, \ldots, x_{k_{m(3)}}]$, where $x_{i_1}, \ldots, x_{k_{m(3)}} \in \mathcal{A}$ and $m(1), m(2), m(3) \geq 2$. Since $\varepsilon(a_\mu) = 0$ with $\mu \in \{1, 2, 3\}$, we get

$$\partial_i([a_1, a_2, a_3]) = \partial_i(a_1)a_2a_3 - \partial_i(a_2)a_1a_3 - \partial_i(a_3)a_1a_2 + \partial_i(a_3)a_2a_1 \tag{4.2}$$

for all $i \in \{1, \ldots, n\}$. We claim that the entries of the $2\times 2$ matrix $\Phi(\partial_i(u_j))$ belong to $\Lambda^{(4)}$ for all $i, j \in \{1, \ldots, n\}$. Indeed, each $u_j \in J_n$ is a $K$-linear combination of Lie commutators of the form $[a_1, a_2, a_3]$, where $a_1, a_2$ and $a_3$ as above, $m(1), m(2), m(3) \geq 2$ and $m(1)+m(2)+m(3) \geq \kappa + 3$. As observed before, for $i_1, \ldots, i_\beta, j_1, \ldots, j_\gamma \in \{1, \ldots, n\}$, with $\beta, \gamma \geq 2$, the entries of $\Phi([x_{i_1}, \ldots, x_{i_\beta}])\Phi([x_{j_1}, \ldots, x_{j_\beta}])$ belong to $\Lambda^{(4)}$. By applying $\Phi$ on (4.2), we get the entries of $\Phi(\partial_i([a_1, a_2, a_3]))$ belong to $\Lambda^{(4)}$ for all $i \in \{1, \ldots, n\}$. Hence the entries of $\Phi(\partial_i(u_j))$ belong to $\Lambda^{(4)}$ for all $u_j \in J_n$, $i, j \in \{1, \ldots, n\}$.

Let $\kappa = 3$. Working modulo $\gamma_7(L_n)$, we have $\psi_3(x_1) \equiv x_1 + v_3 + u_1$ and $\psi_3(x_i) \equiv x_i + u_i$ where each $u_i \in J_n \cap L_n^6$. Since $\psi_3$ is an automorphism of $L_n$, we obtain

$$W = \partial_1(v_3) + (\partial_1(u_1) + \cdots + \partial_n(u_n)) \tag{4.3}$$

is a balanced element. Hence $\mathrm{Tr}(\Phi(W))) = 0$. On the other hand, by (4.2) and since each $u_j \in J_n \cap L_n^6$, we get $\mathrm{Tr}(\Phi(\partial_j(u_j))) \in \Lambda^{(4)}$ for all $j \in \{1, \ldots, n\}$. By applying $\Phi$ on the equation (4.1), we have

$$\Phi(\partial_1(v_3)) = \left(\sum_{\lambda=1}^{3}(-1)^{\lambda-1}\binom{3}{\lambda}\right)\Phi(x_2)[\Phi(x_3), \Phi(x_4)].$$

Since $\sum_{\lambda=1}^{3}(-1)^{\lambda-1}\binom{3}{\lambda} = 1$, we get

$$\begin{aligned}\Phi(\partial_1(v_3)) &= \Phi(x_2)[\Phi(x_3), \Phi(x_4)] \\ &= \Phi(x_2)(\Phi(x_3)\Phi(x_4) - \Phi(x_4)\Phi(x_3)) \\ &= \begin{pmatrix} z^3 & z^3 \\ z^3 & z^3 \end{pmatrix}\end{aligned}$$

and so $\mathrm{Tr}(\Phi(\partial_1(v_3))) = 2z^3$. Therefore, by (4.3), $\mathrm{Tr}(\Phi(W)) \equiv 2z^3$ modulo $\Lambda^{(4)}$ which is a contradiction. Hence $\varphi_3$ is a non-tame automorphism of $C_n$.

Thus we may assume that $\kappa \geq 4$. Let $M_n(A_n)$ be the ring of $n \times n$ matrices over $A_n$. The algebra homomorphism $\Phi$ from $A_n$ into $M_2(K[z])$ defines a homomorphism $\Psi$ from $M_n(A_n)$ into the ring of $2n \times 2n$ partitioned matrices over $K[z]$. Since $\psi_\kappa$ is an automorphism, we get the $2n \times 2n$ partitioned matrix $\Psi((\partial_i(\psi_\kappa(x_j))))$ is invertible over $K[z]$. By applying $\Phi$ on the equation (4.1), we have

$$\Phi(\partial_1(v_\kappa)) = \left(\sum_{\lambda=1}^{\kappa}(-1)^{\lambda-1}\binom{\kappa}{\lambda}\right)\Phi(x_2)[\Phi(x_3), \Phi(x_4)].$$

Since $\sum_{\lambda=1}^{\kappa}(-1)^{\lambda-1}\binom{\kappa}{\lambda} = 1$, we get

$$\begin{aligned}\Phi(\partial_1(\psi_\kappa(x_1)) &= \mathrm{Id}_2 + \Phi(x_2)[\Phi(x_3), \Phi(x_4)] + \Phi(\partial_1(u_1)) \\ &= \mathrm{Id}_2 + \begin{pmatrix} z^3 & z^3 \\ z^3 & z^3 \end{pmatrix} + \Phi(\partial_1(u_1)) \\ &= \begin{pmatrix} 1+z^3 & z^3 \\ z^3 & 1+z^3 \end{pmatrix} + \Phi(\partial_1(u_1))\end{aligned}$$

and, for $i \in \{2, \ldots, n\}$, $\Phi(\partial_i(\psi_\kappa(x_i))) = \mathrm{Id}_2 + \Phi(\partial_i(u_i))$. Since the entries of each $\Phi(\partial_i(u_j))$ belong to $\Lambda^{(4)}$, we have, by working modulo $\Lambda^{(4)}$, $\Psi((\partial_i(\psi_\kappa(x_j))))$ has a non-unit determinant in $K[z]$ which is a contradiction. Therefore, for each $\kappa$, $\varphi_\kappa$ is a non-tame automorphism of $C_n$. □

**Corollary 4.1** *The automorphisms $\omega_1, \omega_2$ and $\omega_3$ of $R_n$ are non-tame.*

*Proof.* Recall that $I_n = \gamma_3(L_n') + (\gamma_3(L_n))'$. Since $I_n \subseteq \gamma_5(L_n)$, it may be deduced from the first part of the above proof that $\omega_1$ and $\omega_2$ are non-tame automorphisms of $R_n$. Throughout this proof, let $v = [x_1, x_2, {}_2x_1, [x_3, x_4]]$. Suppose that $\omega_3$ is induced by an automorphism $\psi_3$ of $L_n$. Then we can write $\psi_3(x_1) = x_1 + v + u_1$ and $\psi_3(x_i) = x_i + u_i$ where each $u_i \in I_n$. Let $b_1 = [x_{i_1}, \ldots, x_{i_{m(1)}}]$ and $b_2 = [x_{j_1}, \ldots, x_{j_{m(2)}}]$, where $x_{i_1}, \ldots, x_{j_{m(2)}} \in \mathcal{A}$ and $m(1), m(2) \geq 3$. Since $\varepsilon(b_\nu) = 0$ with $\nu \in \{1, 2\}$, we get

$$\partial_i([b_1, b_2]) = \partial_i(b_1)b_2 - \partial_i(b_2)b_1 \tag{4.4}$$

for all $i \in \{1, \ldots, n\}$. Working modulo $\gamma_7(L_n)$, we have $\psi_3(x_1) \equiv x_1 + v + u_1$ and $\psi_3(x_i) \equiv x_i + u_i$ where each $u_i \in I_n \cap L_n^6$. Since $I_n \cap L_n^6 = (\gamma_3(L_n') \cap L_n^6) \oplus ((\gamma_3(L_n))' \cap L_n^6)$, we write each $u_i = u_{i1} + u_{i2}$, where $u_{i1} \in \gamma_3(L_n') \cap L_n^6$ and $u_{i2} \in (\gamma_3(L_n))' \cap L_n^6$. Since $\psi_3$ is an automorphism of $L_n$, we obtain

$$U = \partial_1(v) + (\partial_1(u_{11}) + \cdots + \partial_n(u_{n1})) + (\partial_1(u_{12}) + \cdots + \partial_n(u_{n2}))$$

is a balanced element. As in the proof of Proposition 4.1, $\mathrm{Tr}(\Phi(U)) = 0$. On the other hand, by (4.2), (4.4) and since each $u_i \in I_n \cap L_n^6$, we get $\mathrm{Tr}(\Phi(\partial_i(u_i))) \in \Lambda^{(4)}$ for all $i \in \{1, \ldots, n\}$. As in the proof of Proposition 4.1, $\mathrm{Tr}(\phi(\partial_1(v))) = 2z^3$. Therefore $\mathrm{Tr}(\Phi(U)) \equiv 2z^3$ modulo $\Lambda^{(4)}$ which is a contradiction. Hence $\omega_3$ is a non-tame automorphism of $R_n$. □